\documentclass{article}
\usepackage{amssymb}
\usepackage{amsfonts}
\usepackage{graphicx}
\usepackage{amsmath}

\newtheorem{theorem}{Theorem}

\newtheorem{definition}[theorem]{Definition}

\newtheorem{proposition}[theorem]{Proposition}
\newtheorem{remark}[theorem]{Remark}

\input{tcilatex}
\begin{document}

\title{Riemannian geometry as a \\
curved pre-homogeneous geometry }
\author{ Erc\"{u}ment Orta\c{c}gil}
\maketitle

\begin{abstract}
We define a Riemannian structure as a pre-homogeneous geometric structure
with curvature $\mathcal{R}.$ We show that $\mathcal{R}=0$ if and only if
the underlying metric has constant curvature. We define pre-homogeneous
geometric structures and pose some problems.
\end{abstract}

\section{Introduction}

This note is the continuation of [24], [3], [2] and its main purpose is to
carry out the program outlined in the introduction of [2] in the case of
Riemannian and affine structures. So we will start by recalling this program
in more technical detail than [2].

Let $G$ be a connected Lie transformation group which acts effectively and
transitively on a connected and smooth manifold $M.$ We call this data a
Klein geometry and denote it by $(G,M).$We fix some $p,q\in M,$ $g\in G$
with $g(p)=q,$ and define $\mathcal{H}_{k}\overset{def}{=}\{h\in G\mid
j_{k}(h)^{p,q}=j_{k}(g)^{p,q},$ $k\geq 0\}$ where $j_{k}(g)^{p,q}$ denotes
the $k$-jet of the transformation $g$ with source at $p$ and target at $q.$
We call $j_{k}(g)^{p,q}$ a $k$-arrow (from $p$ to $q)$. So $\mathcal{H}_{k}$
consists of all $h\in G$ such that the $k$-arrow of $h$ coincides with the $%
k $-arrow of $g$ $($for $k=0,$ this simply means $h(p)=g(p)).$ Clearly $%
\{g\}\subset \mathcal{H}_{k+1}\subset \mathcal{H}_{k}.$ In [3], we showed
that this sequence stabilizes at $\{g\}$ and the smallest integer $m$ with $%
\mathcal{H}_{m}=\{g\}$ does not depend on $p,q$ and $g.$ Therefore, any $%
g\in G$ is determined by its $m$-arrow $j_{m}(g)^{p,q}$ for any $p,q$ and
this holds in particular for the elements of the stabilizer $H_{p}\overset{%
def}{=}\{g\in G\mid g(p)=p\}.$ The integer $m$ is called the geometric order
of the Klein geometry $(G,M)$ in [3]$.$ Clearly, $m=0$ if and only if $G$
acts simply transitively. We showed in [3] that $m$ can be arbitrarily large
even if $M$ is compact and actions whose stabilizers are parabolic subgroups
are the prototypes of this situation.

There is no curvature in the above picture in the sense that any $m$-arrow $%
j_{m}(g)^{p,q}$ integrates uniquely to a local (in fact global)
diffeomorphism which is the transformation $g$. Shortly, the groupoid of $m$%
-arrows integrates to a pseudogroup (in fact to a Lie group). So the
question arises how we should \textquotedblleft curve\textquotedblright\ the
above picture, that is, what is curvature? However we look at curvature, it
must deform the \textquotedblleft symmetric\textquotedblright\ object $(G,M)$
into a \textquotedblleft lumpy\textquotedblright\ object. All the existing
approaches to the concept of curvature (see, for instance, [14], [29], [6],
[8]) circle around this fundamental idea and this note and [2] are no
exceptions. It is commonly accepted today (at least in Riemannian geometry,
see [15] and the recent work [8] for parabolic geometries) that this lumpy
object is a principal bundle $P\rightarrow N$ with structure group $H\simeq
H_{p},$ $\dim N=\dim M,$ together with some extra structure on $P\rightarrow
N$ , like a torsionfree connection. A different approach is taken in [6]
which is more in the spirit of this note and [2]. To explain our lumpy
object $\mathcal{G}(N)$, we continue with the setting of the above paragraph
and let $\mathcal{U}^{p,p}(r)$ denote the group consisting of $r$-arrows of 
\textit{all }local diffeomorphisms with source and target at $p\in M,$ $%
r\geq 0.$ We fix a coordinate system around $p$ which identifies $\mathcal{U}%
^{p,p}(r)$ with the jet group $G_{r}(n)$ where $n=\dim M$. Since $%
j_{m}(h)^{p,p}$ determines $h\in H_{p}$ (for a technical reason to be
explained in Remark 1 and Section 7, we assume here that $M$ is simply
connected), $H_{p}$ injects onto a subgroup of $G_{m}(n)$ which we denote
also by $H_{p}$. Since $h\in H_{p}$ determines $j_{m+1}(h)^{p,p},$ we also
have an injection of $H_{p}$ onto a subgroup $G_{m+1}(n)$ which we denote by 
$\varepsilon H_{p}.$ So we have the commutative diagram

\begin{equation}
\begin{array}{ccc}
G_{m+1}(n) & \overset{\pi }{\longrightarrow } & G_{m}(n) \\ 
\cup &  & \cup \\ 
\varepsilon H_{p} & \overset{\pi }{\longrightarrow } & H_{p}%
\end{array}%
\end{equation}%
where $\pi $ is the projection homomorphism induced by the projection of
jets and the vertical inclusions depend on the coordinate system around $p$
which we fixed. The restriction of $\pi $ to $\varepsilon H_{p}$ is a
bijection in (1) and $\varepsilon =(\pi _{\mid H_{p}})^{-1}$ so that $\pi
\circ \varepsilon =id$. It is easy to check that the conjugacy class of $%
\varepsilon H_{p}$ inside $G_{m+1}(n)$ is independent of the coordinate
system around $p$ and therefore also independent of $p$ since any two
stabilizers are conjugate. We denote this conjugacy class by $\{G,M,H\}$
where $H$ denotes some \textquotedblleft easy\textquotedblright\
representative of the isomorphism class of the stabilizers $H_{p},$ $p\in M$%
. For instance in Eucledean geometry, $H=O(n)$, $G=O(n)\ltimes \mathbb{R}%
^{n} $ and $M=\mathbb{R}^{n}.$ We call $\{G,M,H\}$ the vertex connection of
the Klein geometry $(G,M).$ In disguise, the vertex connection $\{G,M,H\}$
is a $PDE$ for $\dim M\geq $ $2$ and $ODE$ for $\dim M=1$ and plays a
fundamental role in the theory. Two Klein geometries $(G_{1},M_{1}),$ $%
(G_{2},M_{2})$ with $\dim M_{1}=\dim M_{2},$ $\dim H_{1}=\dim H_{2}$ and $%
m_{1}=m_{2}$ may define the same vertex connection $\{G,M,H\}$ but $G_{1}$
and $G_{2}$ may have nonisomorphic Lie algebras, that is, even though $%
\{G,M,H\}$ by definition determines the isomorphism class of $H$, it does
not completely determine $G$ locally. In Riemannian geometry, there are only
three and in affine geometry only one such Lie algebra up to isomorphism.
The other extreme case occurs in [2] for $m=0$ where \textit{all }Lie
algebras arise.

Now, in search for a geometric structure $\mathcal{G}(N)$ which is the
curved analog of $(G,M),$ the main idea of the reference [22] in [2] is to
dispense with $(G,M)$ but retain the vertex connection $\{G,M,H\},$ start
with some smooth manifold $N$ with $\dim N=\dim M$, and try to reconstruct
the action of $G$ on $M$ from a transitive Lie groupoid $\mathcal{G}(N)$ of
order $m+1$ on $N$ with the property that the conjugacy class of the vertex
groups of $\mathcal{G}(N)$ is $\{G,M,H\}$ (see Definition 14). As required, $%
(m+1)$-arrows of $\mathcal{G}(N)$ integrate to a pseudogroup if and only if
the curvature $\mathcal{R}=0$ so that $\mathcal{R}$ is the obstruction to
construct the action of $G$ locally. The lift of this pseudogroup globalizes
to a Lie group $\widetilde{G}$ on the universal covering space $\widetilde{N%
\text{ }}$of $N$ and the Klein geometry $(\widetilde{G},\widetilde{N})$
defines the vertex connection that we started with. In this way we obtain
the well known uniformization theorems for Riemannian and affine structures
(Propositions 9, 10). The remarkable fact is that $\mathcal{R}$ and the
curvature in the formalism of connections on principal bundles are different
objects, even in Riemannian geometry, as can be seen from our abstract. We
hope that these differences will become more transparent in this note.
Another surprising fact is that we mention torsion, the Levi-Civita
connection and covariant differentiation a few times in this note only for
convinience.

The above idea is worked out in detail for $m=0$ in [2]. Even in this
simplest case (or the most complicated case depending on our view), this
approach gives rise to some new possibilities and questions with wide scope
and subtlety. The purpose of this note is to show how this program works out
for $m=1$ in the case of Riemannian and affine structures and indicate how
it generalizes in a straightforward way to all pre-homogeneous structures of
arbitrary geometric order $m$.

This note is organized as follows. In Section 2 we fix some
\textquotedblleft easy\textquotedblright\ representatives of the Riemannian
and affine vertex connections and derive some elementary formulas. In
Section 3 we define a Riemannian structure as a subgroupoid $\mathcal{G}%
_{2}\subset \mathcal{U}_{2}$ with algebroid $\mathfrak{G}_{2}\subset J^{2}T$
in accordance with the general theory and using the formulas in Section 2,
we express these submanifolds locally as the zero set of some functions. We
define the complete integrability of $\mathcal{G}_{2}$ and $\mathfrak{G}%
_{2}, $ a concept which plays a fundamental role in this note. In Section 4
we outline a proof that the Lie's third theorem is equivalent to the well
known uniformization theorem in Riemannian geometry (Propositions 9, 10) .
In Section 5 we define the algebroid and groupoid curvatures $\mathfrak{R}$
and $\mathcal{R}$ show that their vanishing is equivalent to complete
integrability and also to the constant curvature condition of the metric
(Propositions 11, 12). In Section 6 we take a brief look at affine
structures. Since all pre-homogeneous structures are studied on the same
footing, all our propositions in Sections 4, 5, 6 are identical with those
in [2], except that we do not touch characteristic classes here. It has
become apparent to us that we will never be able to finish the reference
[22] in [2], so we decided to give here the definition of a pre-homogeneous
structure and formulate some problems which we believe are fundamental for
the theory, which is the content of Section 7.

\section{Riemannian and affine vertex connections}

Let $G_{k}(n)$ denote the $k$'th order jet group in $n$-variables. The
elements of $G_{k}(n)$ are $k$-jets $j_{k}(f)^{o,o}$ of local
diffeomorphisms $f$ of $\mathbb{R}^{n}$ with $f(o)=o$ where $o$ is the
origin. The composition $\circ $ is defined by $j_{k}(f)^{o,o}\circ
j_{k}(g)^{o,o}\overset{def}{=}j_{k}(f\circ g)^{o,o}$. The projection $\pi $
of jets gives the exact sequence of Lie groups 
\begin{equation}
0\longrightarrow K_{k+1,k}(n)\longrightarrow G_{k+1}(n)\overset{\pi }{%
\longrightarrow }G_{k}(n)\longrightarrow 1
\end{equation}%
where $K_{k+1,k}(n)$ is a vector group. In this note, all projections
induced by the projection of jets will be denoted by $\pi $ and all
splittings by $\varepsilon .$ Up to Section 7, we will have $k=1$ in (2). We
refer to [30] for some basic structure theorems for $G_{k}(n)$ and to [22]
for an explicit matrix representation of $G_{k}(n)$.

In the coordinates of $\mathbb{R}^{n}$, an element $a\in G_{2}(n)$ is of the
form $(a_{j}^{i},a_{jk}^{i})$ and the chain rule of differentiation shows
that the group operation is given by

\begin{equation}
(a_{j}^{i},a_{jk}^{i})(b_{j}^{i},b_{jk}^{i})=(a_{s}^{i}b_{j}^{s},a_{s}^{i}b_{jk}^{s}+a_{st}^{i}b_{j}^{s}b_{k}^{t})
\end{equation}%
We use summation convention in (3). For simplicity of notation, we denote $%
a=(a_{j}^{i},a_{jk}^{i})$ by $(a_{1},a_{2})$ and write the group operation
(3) as 
\begin{equation}
ab=(a_{1},a_{2})(b_{1},b_{2})=(a_{1}b_{1},a_{1}b_{2}+a_{2}(b_{1})^{2})
\end{equation}

Clearly, $a_{j}^{i}=\delta _{j}^{i},$ $a_{jk}^{i}=0$ defines the identity
which we denote by $(I,0).$ It is easy to check that $%
(a_{j}^{i},a_{jk}^{i})^{-1}=((a^{-1})_{j}^{i},$ $%
-(a^{-1})_{s}^{i}a_{tr}^{s}(a^{-1})_{j}^{t}(a^{-1})_{k}^{r})$ which we write
as

\begin{equation}
(a_{1},a_{2})^{-1}=(a_{1}^{-1},-a_{1}^{-1}a_{2}(a_{1}^{-1})^{2})
\end{equation}%
using our notation in (4). We have $(I,a)(I,b)=(I,a+b)$ and $%
(I,a)^{-1}=(I,-a).$ We can write (2) now as $(I,a_{2})\longrightarrow
(a_{1},a_{2})\longrightarrow a_{1}.$

We now define a splitting $\varepsilon :O(n)\rightarrow G_{2}(n).$ Let $G(0)$
denote the isometry group of $\mathbb{R}^{n},$ that is, $G(0)=O(n\mathbb{%
)\ltimes R}^{n}$ where $\ltimes $ denotes semidirect product. For $g=(\xi
,a)\in G(0),$ $x\in \mathbb{R}^{n},$ we have $(gx)^{i}=\xi
_{s}^{i}x^{s}+b^{i}.$ Therefore, if $h(p)=p$ for some $h\in G(0)$ and $p\in 
\mathbb{R}^{n}$, then $\left[ j_{1}(h)^{p,p}\right] _{j}^{i}=\xi _{j}^{i}$
and $\left[ j_{2}(h)^{p,p}\right] _{j}^{i}=0$. Therefore any $h\in G(0)$
which stabilizes $p$ (therefore any $g\in G(0))$ is determined by its $1$%
-arrow $j_{1}(g)^{p,p}$ (note the crucial role of translations!) and it
follows that $m=1$ where $m$ is the geometric order of the Klein geometry $%
(G(0),O(n\mathbb{)}).$

\begin{remark}
In the definition of the geometric order of the Klein geometry $(G,M)$ in
[3], the connectedness of $G$ is used to ensure that the adjoint
representation of $H_{p}$ on the Lie algebra of $G$ is faithful (see Lemma
5.1 in [3]). Henceforth we will always assume this latter condition so that
geometric order is defined. Also, if $M$ is simply connected, then the
geometric order of $(G,M)$ is equal to the infinitesimal order of $\mathfrak{%
(g,h)}$ where $\mathfrak{h}$ is the Lie algebra of some stabilizer (see [3]
for the geometric and infinitesimal orders and the second paragraph of
Section 7 of this note)
\end{remark}

Thus $O(n)$ injects into $G_{2}(n)$ as%
\begin{eqnarray}
\varepsilon &:&O(n)\rightarrow G_{2}(n) \\
&:&a\rightarrow (a,0)  \notag
\end{eqnarray}

Since the above derivation does not use the orthogonality of the matrix $\xi 
$, we may replace $O(n\mathbb{)}$ with $G_{1}(n)$ and $\varepsilon $ splits
also $G_{1}(n)$ (in fact, any subgroup of $G_{1}(n))$ by the same formula in
(6)). Therefore $m=1$ also for the affine group $A=G_{1}(n)\ltimes \mathbb{R}%
^{n}.$

Now we denote the conjugacy classes of $\varepsilon O(n)$ and $\varepsilon
G_{1}(n)$ inside $G_{2}(n)$ by $\{G(0),\mathbb{R}^{n},SO(n\mathbb{)\}}$ and $%
\{A,\mathbb{R}^{n},G_{1}(n)\}$ respectively. It is crucial to observe how
the independence of $p$ and the coordinates around $p$ explained in the
Introduction is \textquotedblleft trivialized\textquotedblright\ by
translations which will not be at our disposal if we replace $\mathbb{R}^{n}$
with some arbitrary $M.$

These two objects play a fundamental role in this note, so we make

\begin{definition}
The conjugacy classes $\{G(0),\mathbb{R}^{n},O(n\mathbb{)\}}$ and $\{A,%
\mathbb{R}^{n},G_{1}(n)\}$ are the vertex connections of the Klein
geometries $(G(0),\mathbb{R}^{n})$ and $(A,\mathbb{R}^{n})$ respectively.
\end{definition}

We denote the vertex connections in Definition 2 by $\mathbf{R}$ and $%
\mathbf{A}$ and call them Riemannian and affine respectively. So $%
\varepsilon O(n)$ and $\varepsilon G_{1}(n)$ are some \textquotedblleft
easy\textquotedblright\ representatives of $\mathbf{R}$ and $\mathbf{A.}$
These choices are irrelevant from a theoretical standpoint but greatly
simplify local computations.

Now (4) gives $%
(a_{1},a_{2})=(a_{1},0)(a_{1},0)^{-1}(a_{1},a_{2})=(a_{1},0)(I,a_{1}{}^{-1}a_{2}) 
$ which expresses the semidirect product structure

\begin{equation}
G_{2}(n)=G_{1}(n)\ltimes K_{2,1}(n)
\end{equation}

In fact, $G_{1}(n)$ splits inside $G_{k}(n)$ for all $k\geq 1,$ $n\geq 1$ in
the same way, that is, $G_{k}(n)=G_{1}(n)\ltimes K_{k,1}$. In more abstract
terms, an algebraic group is the semidirect product of its maximal reductive
and maximal nilpotent subgroups and the decomposition $G_{k}(n)=G_{1}(n)%
\ltimes K_{k,1}$ is a special case ([30], Theorem 2.6). We believe that the
splitting of (2) for all $n\geq 1$ occurs only for $k=1$ and (2) never
splits for $n\geq 2,$ $k\geq 2.$ We will see in Section 7 that the
Schwarzian derivative arises from the splitting of (2) for $n=1,$ $k=2.$ All
splittings of (2) for $n=1$ and arbitrary $k$ are determined in [27] on the
level of Lie algebras. We also refer to [30] for a detailed study of the
solvable Lie group $G_{k}(1)$. Now, even though (2) splits very rarely,
there exist an abundance of splittings \textit{inside }(2) for arbitrarily
large values of $n$, $k$ arising from Klein geometries $(G,M)$ where $n=\dim
M,$ $k=$ geometric order of $(G,M).$ These splittings form the backbone of
the present approach as explained in the Introduction.

Now let $G(1)\overset{def}{=}O(n+1)$, $G(-1)\overset{def}{=}O(1,n)$ so that
we have

\begin{equation}
G(1)/O(n)=\mathbb{S}^{n},\text{ \ \ \ }G(0)/O(n)=\mathbb{R}^{n},\text{ \ \ \ 
}G(-1)/O(n)=\mathbb{H}^{n}
\end{equation}

The geometric order $m=1$ for the Klein geometries in (8). More generally,
if $H\neq \{e\}$ is compact and $G/H$ is simply connected, then $m=1.$ This
is equivalent to the statement that the isotropy representation of $H$ is
faithful. It follows that the vertex connections $\{G(1),\mathbb{S}%
^{n},O(n)\}$, $\{G(-1),\mathbb{H}^{n},O(n))$ are defined as explained in the
Introduction. We claim

\begin{equation}
\{G(1),\mathbb{S}^{n},O(n)\}=\{G(0),\mathbb{R}^{n},O(n\mathbb{)\}=}\{G(-1),%
\mathbb{H}^{n},O(n)\}
\end{equation}

We will see in Section 4 that (9) is a consequence of the uniformization
theorem in Riemannian geometry.

We now take a closer look at $\{A,\mathbb{R}^{n},G_{1}(n)\}.$

\begin{proposition}
Let $\sigma _{1},\sigma _{2}$ $:G_{1}(n)\rightarrow G_{2}(n)$ be two group
homomorphisms satisfying $\pi \circ \sigma _{i}=id.$ Then $\sigma
_{1}G_{1}(n)=k(\sigma _{2}G_{1}(n))k^{-1}$ for some $k\in
K_{2,1}(n)\vartriangleleft G_{2}(n).$
\end{proposition}

To prove the assertion, let $\sigma :G_{1}(n)\rightarrow G_{2}(n)$ be any
such homomorphism. So $\sigma (a)=(a,\phi (a))$ for some function $\phi .$
Now $\sigma (ab)=\sigma (a)\sigma (b)$ and (4) give

\begin{equation}
\phi (ab)=a\phi (b)+\phi (a)b^{2}\text{ \ \ }a,b\in G_{1}(n)
\end{equation}%
We successively let $a=\lambda I$ and $b=\lambda I$ in (10) and get

\begin{equation}
\phi (\lambda b)=\lambda \phi (b)+\phi (\lambda )b^{2}\text{ \ \ },\text{\ \
\ }\phi (b\lambda )=b\phi (\lambda )+\phi (b)\lambda ^{2}\text{ \ \ \ }b\in
G_{1}(n)
\end{equation}%
where $\lambda $ denotes $\lambda I$ in (11). Since $\lambda b=b\lambda ,$
(11) gives $\lambda \phi (b)+\phi (\lambda )b^{2}=b\phi (\lambda )+\phi
(b)\lambda ^{2}$. Solving for $\phi (b),$ we obtain

\begin{eqnarray}
(b,\phi (b)) &=&\left( b,\frac{1}{\lambda ^{2}-\lambda }(\phi (\lambda
I)b^{2}-\text{\ }b\phi (\lambda I)\right) \text{\ \ \ \ }\lambda \in \mathbb{%
R},\text{ }\lambda \neq 0,1,\text{ \ }b\in G_{1}(n)  \notag \\
&=&(I,\frac{\phi (\lambda I)}{\lambda ^{2}-\lambda })(b,0)(I,\frac{\phi
(\lambda I)}{\lambda ^{2}-\lambda })^{-1}\text{ \ \ \ }
\end{eqnarray}

We observe that the $RHS$ of (12) is independent of $\lambda $. Setting $k=%
\frac{\phi (\lambda I)}{\lambda ^{2}-\lambda }$, (12) becomes $\sigma
(b)=(b,\phi (b))=(I,k)\varepsilon (b)(I,-k)$. Therefore any splitting $%
\sigma $ is conjugate to $\varepsilon $, which proves the statement.

Note that the above proof works if we replace $G_{1}(n)$ by any subgroup $%
L\subset G_{1}(n)$ as long as $L$ contains some $\lambda I$ with $\lambda
\neq 1$. Therefore it works for $O(n)$ and $O(1,n)$ since $-I\in O(n),$ $%
O(1,n).$ Also, note that we do not even assume the continuity of $\sigma $
but deduce the very strong conclusion that $\sigma (b)$ is a quadratic
polynomial in $b.$

Now our purpose is to express the right cosets of $\varepsilon (O(n))$ in $%
G_{2}(n)$ as the zero set of some functions. We have%
\begin{eqnarray}
(a_{1},a_{2})(b_{1},b_{2})^{-1} &\in &\varepsilon (O(n))\Longleftrightarrow
(a_{1},a_{2})(b_{1}^{-1},-b_{1}^{-1}b_{2}(b_{1}^{-1})^{2}\in \varepsilon
(O(n))  \notag \\
&\Longleftrightarrow
&(a_{1}b_{1}^{-1},-a_{1}b_{1}^{-1}b_{2}(b_{1}^{-1})^{2}+a_{2}(b_{1}^{-1})^{2}\in \varepsilon (O(n))
\notag \\
&\Longleftrightarrow &(a_{1}b_{1}^{-1})^{T}(a_{1}b_{1}^{-1})=I,\text{ \ \ }%
-a_{1}b_{1}^{-1}b_{2}(b_{1}^{-1})^{2}+a_{2}(b_{1}^{-1})^{2}=0  \notag \\
&\Longleftrightarrow &(b_{1}^{T})^{-1}a_{1}^{T}a_{1}b_{1}^{-1}=I,\text{ \ \ }%
a_{1}b_{1}^{-1}b_{2}=a_{2}  \notag \\
&\Longleftrightarrow &a_{1}^{T}a_{1}=b_{1}^{T}b_{1},\text{ \ \ }%
a_{1}^{-1}a_{2}=b_{1}^{-1}b_{2}
\end{eqnarray}

We define the functions $F_{1},$ $F_{2}$ by $F_{1}(a_{1},a_{2})\overset{def}{%
=}a_{1}^{T}a_{1},$ $F_{2}(a_{1},a_{2})\overset{def}{=}a_{1}^{-1}a_{2}$, set $%
F\overset{def}{=}(F_{1},F_{2})$ and rewrite (13) as

\begin{equation}
ab^{-1}\in \varepsilon (O(n))\Longleftrightarrow F(a)=F(b)
\end{equation}%
In more detail, $F_{1}$ has components $F_{jk}:G_{2}(n)\rightarrow \mathbb{R}%
,$ $1\leq j,k\leq n,$ defined by

\begin{equation}
F_{jk}(a_{1},a_{2})=a_{j}^{s}a_{k}^{s}=(a_{1}^{T}a_{1})_{jk}
\end{equation}%
where $s$ is summed in (15), and $F_{2}$ has components $F_{jk}^{i}:G_{2}(n)%
\rightarrow \mathbb{R},$ $1\leq i,j,k\leq n,$ defined by

\begin{equation}
F_{jk}^{i}(a_{1},a_{2})=(a^{-1})_{s}^{i}a_{jk}^{s}=(a_{1}^{-1}a_{2})_{jk}^{i}
\end{equation}

So $F:G_{2}(n)\rightarrow \mathbb{R}^{s},$ $s=\dim G_{2}(n)$ and $F$ has
constant rank $r=\dim G_{2}(n)-\dim \varepsilon (O(n))$. Thus the surjective
map $F:G_{2}(n)\rightarrow \mathbb{R}^{r}$ has right cosets of $\varepsilon
(O(n))$ as fibers. If we replace $\varepsilon (O(n))$ by $\varepsilon
(G_{1}(n))$, the function $F$ will be defined only by $F_{2}$ as the
condition imposed by $F_{1}$ will be redundant. We will continue to use the
notation $F$ also in this case. This point will be relevant in Section 6.

Recall that $G_{2}(n)$ is an algebraic group, $\varepsilon (O(n)),$ $%
\varepsilon (G_{1}(n))\subset G_{2}(n)$ are algebraic subgroups and $F$ is a
polynomial. On the other hand, if $o$ denotes the coset of $\varepsilon
(O(n) $ (or $\varepsilon (G_{1}(n))$, we can always express the cosets 
\textit{near }$o$ as the zero set of some smooth functions.

Now suppose $ab=c$ in $G_{2}(n).$ We have $%
F_{1}(c)=(a_{1}b_{1})^{T}a_{1}b_{1}=b_{1}^{T}(a_{1}^{T}a_{1})b_{1}=b_{1}^{T}F_{1}(a)b_{1}, 
$ that is

\begin{equation}
c_{jk}=(a_{1}^{T}a_{1})_{st}b_{j}^{s}b_{k}^{t}
\end{equation}

Similarly, $%
F_{2}(c)=F_{2}(ab)=(a_{1}b_{1})^{-1}(a_{1}b_{2}+a_{2}b_{1})=b_{1}^{-1}b_{2}+b_{1}^{-1}(a_{1}^{-1}a_{2})b_{1}=F_{2}(b)+b_{1}^{-1}F_{2}(a)b_{1}, 
$ that is,

\begin{equation}
c_{jk}^{i}=a_{tr}^{s}b_{j}^{t}b_{k}^{r}(b^{-1})_{s}^{i}+(b^{-1})_{s}^{i}b_{jk}^{s}
\end{equation}

Finally, it is easy to check that (17), (18) satisfy the group law, that is,
the composition $F(a)\rightarrow F(ab)\rightarrow F((ab)c)$ is the same as $%
F(a)\rightarrow F(a(bc)).$

\section{Riemannian structures}

Let $M$ be a smooth and connected manifold with $\dim M=n.$ Let $p,q\in M$
and $\mathcal{U}_{k}^{p,q}$ denote the set of \textit{all }$k$-jets $%
j_{k}(f)^{p,q}$ of local diffeomorphisms with source at $p$ and target at $%
q. $ We call $j_{k}(f)^{p,q}$ a $k$-arrow (from $p$ to $q).$ The composition
of local diffeomorphisms induces a composition $\mathcal{U}_{k}^{q,r}\times $
$\mathcal{U}_{k}^{p,q}\rightarrow $ $\mathcal{U}_{k}^{p,r}.$ We define the
set $\mathcal{U}_{k}\overset{def}{=}\cup _{p,q\in M}\mathcal{U}_{k}^{p,q}$.
The smooth structure of $M$ induces a natural smooth structure on $\mathcal{U%
}_{k}$ as follows. For two coordinate patches $(U,x^{i})$, $(V,y^{i})$ on $M$%
, any $k$-arrow $f^{p,q}\in \mathcal{U}_{k}$ with $p\in U$ and $q\in V$ has
the unique representation $%
(x^{i},y^{i},y_{j_{1}}^{i},y_{j_{2}j_{1}}^{i},...,y_{j_{k}...j_{1}}^{i})$
where $x^{i}$ and $y^{i}$ are the coordinates of $p,q$ respectively$.$ With
this differentiable structure, $\mathcal{U}_{k}$ \ becomes a transitive Lie
equation in finite form which is a very special groupoid (see [28], [25] for
Lie equations in finite and infinitesimal forms and [18] and the references
therein for general Lie groupoids and algebroids). We call $\mathcal{U}_{k}$
the universal groupoid on $M$ of order $k.$ Since a $0$-arrow is an ordered
pair $(p,q)$, $\mathcal{U}_{0}$ is the pair groupoid $M\times M.$ Note that
a choice of coordinates around $p\in M$ identifies the vertex group $%
\mathcal{U}_{k}^{p,p}$with the jet group $G_{k}(n)$ for $k\geq 1$ and a
change of these coordinates conjugates this identification with the $k$-jet
of the coordinate change at $p.$ The projection of jets induces a projection 
$\pi :\mathcal{U}_{k+1}\rightarrow \mathcal{U}_{k}$ and $\pi $ is a morphism
of groupoids, that is, it preserves the composition and inversion of arrows.

Now let $\mathcal{G}_{k}\subset \mathcal{U}_{k}$ be a transitive
subgroupoid. This means that the set $\mathcal{G}_{k}^{p,q}$ of $k$-arrows
of $\mathcal{G}_{k}$ is nonempty for all $p,q\in M,$ the $k$-arrows of $%
\mathcal{G}_{k}$ are closed under composition and inversion, and $\mathcal{G}%
_{k}\subset \mathcal{U}_{k}$ is an imbedded submanifold (see [18] for
details). The Lie subgroup $\mathcal{G}_{k}^{p,p}\subset \mathcal{U}%
_{k}^{p,p}$ is called the vertex group of $\mathcal{G}_{k}$ at $p.$

We now fix some $p\in M$ and choose some coordinates $(x^{i})$ around $p.$
This choice identifies $\mathcal{U}_{k}^{p,p}$ with $G_{k}(n)$ and therefore
identifies $\mathcal{G}_{k}^{p,p}$ with a subgroup $i_{x}(\mathcal{G}%
_{k}^{p,p})\subset $ $G_{k}(n).$ A change of coordinates $%
f:(x^{i})\rightarrow (y^{i})$ around $p$ conjugates $i_{x}(\mathcal{G}%
_{k}^{p,p})$ with the $k$-jet of $f$ at $p.$ Therefore $i_{x}(\mathcal{G}%
_{k}^{p,p})$ and $i_{y}(\mathcal{G}_{k}^{p,p})$ belong to the same conjugacy
class in $G_{k}(n).$ This conjugacy class does not depend also on the choice
of the point $p$. To see this, let $q\in M$ any other point, choose \textit{%
any }$k$-arrow $f^{p,q}\in $ $\mathcal{U}_{k}^{p,q}$ and fix \textit{any }%
coordinate system $(U,x^{i})$ around $p.$ Now there exists a local
diffeomorphism $g$ with $g(p)=q$ such that $j_{k}(g)^{p,q}=f^{p,q}$ by the
definition of $f^{p,q}.$ So $g$ defines a coordinate system $(V,y^{i})$
around $q$ and the representation of $f^{p,q}$ with respect to $(U,x^{i}),$ $%
(V,y^{i})$ is $(\overline{x}^{i},\overline{x}^{i},\delta _{j}^{i})$ where $%
p=(\overline{x}^{i}).$ Therefore, $(V,y^{i})$ imbeds $\mathcal{G}_{k}^{q,q}$
into $G_{k}(n)$ in the same way as $(U,x^{i})$ imbeds $\mathcal{G}%
_{k}^{q,q}. $ We denote this conjugacy class by $\{\mathcal{G}_{k}\}$ which
is a common property of all the vertex groups of $\mathcal{G}_{k}.$

Recalling the definition of the Riemannian vertex connection $\mathbf{R}$ in
Section 2, we now make the following

\begin{definition}
A Riemannian structure on $M$ is a transitive subgroupoid $\mathcal{G}%
_{2}\subset \mathcal{U}_{2}$ such that $\{\mathcal{G}_{2}\}=\mathbf{R}$
\end{definition}

Note that a Riemannian structure is a second order structure according to
Definition 4 in the same way as a parallelizable manifold is a first order
structure in [2]. We define $\mathcal{G}_{1}\overset{def}{=}\pi \mathcal{G}%
_{2}$ and call $\mathcal{G}_{1}$ the underlying metric structure. We may
have $\mathcal{G}_{2}\neq \overline{\mathcal{G}_{2}}$ but $\pi \mathcal{G}%
_{2}=\pi \overline{\mathcal{G}_{2}}$ that is, two different Riemannian
structures may have the same underlying metric. By the definitions of $%
\mathcal{G}_{2},$ $\{\mathcal{G}_{2}\},$ $\mathbf{R,}$ for any $p\in M,$
there exists a coordinate system $(U,x^{i})$ around $p$ such that the vertex
group $\mathcal{G}_{2}^{p,p}$ imbeds into $G_{2}(n)$ as (6). We call $%
(U,x^{i})$ regular at $p.$ It will become clear below that regular and
geodesic coordinates agree to the first order, but not neccessarilly to the
second order.

Now, the projection $\pi :\mathcal{G}_{2}\rightarrow \mathcal{G}_{1}$
induces an isomorphism on the vertex groups by the definition of $\mathcal{G}%
_{2}.$ This fact implies that $\pi $ is an isomorphism of groupoids, that
is, above any $1$-arrow of $\mathcal{G}_{1}$, there is a unique $2$-arrow of 
$\mathcal{G}_{2}.$ Indeed, $\pi (f^{p,q})=\pi (g^{p,q})$ if and only if $\pi
\left( (f^{p,q})^{-1}\circ g^{p,q}\right) =j_{1}(id)^{p,p},$ but $%
j_{2}(id)^{p,p}$ is the unique point satisfying $\pi \left(
j_{2}(id)^{p,p}\right) =j_{1}(id)^{p,p},$ so we conclude $f^{p,q}=g^{p,q}.$
Thus we have the the commutative diagram

\begin{equation}
\begin{array}{ccc}
\mathcal{U}_{2} & \overset{\pi }{\longrightarrow } & \mathcal{U}_{1} \\ 
\cup &  & \cup \\ 
\mathcal{G}_{2} & \overset{\pi }{\longrightarrow } & \mathcal{G}_{1}%
\end{array}%
\end{equation}%
where $\pi _{\mid \mathcal{G}_{2}}$ is an isomorphism with inverse $%
\varepsilon \overset{def}{=}(\pi _{\mid \mathcal{G}_{2}})^{-1}$ and as
remarked above, the restriction of $\varepsilon $ to the vertex groups
\textquotedblleft looks like\textquotedblright\ (6) in regular coordinates.
It is crucial to observe that a groupoid is by no means determined by its
vertex groups. An extreme case occurs in [2] where the vertex groups are
trivial, but clearly $M$ can be parallelized in quite different ways.
However, we will see in Sections 4, 5, 6 that the vertex groups severely
restrict Riemannian and affine structures when curvature vanishes.

The above local coordinates on $\mathcal{U}_{2}$ show that $\dim \mathcal{U}%
_{2}=2\dim M+\dim G_{2}(n)$ and $\dim \mathcal{G}_{2}=2\dim M+\dim O(n).$
Now our purpose is to express the submanifold $\mathcal{G}_{2}\subset 
\mathcal{U}_{2}$ locally as the zero set of some functions so that we can
take a closer look at $\mathcal{G}_{2}$ and make some explicit local
computations. So we fix some base point $e\in M$ and some regular
coordinates around $e$ once and for all. Changing these choices will
conjugate our formulas by some \textquotedblleft constant\textquotedblright\ 
$2$-arrow which does not depend on the base variables, so that our formulas
will remain essentially the same when we differentiate them (the reader may
keep track of this in what follows). Now let $(U,x^{i})$ be some arbitrary
coordinate patch on $M.$ For each $p\in U,$ we choose a $2$-arrow of $%
\mathcal{G}_{2}$ with source at $p$ and target at $e.$ This local (smooth)
section $s$ has the coordinate representation $%
s(x)=(x^{i},e^{i},s_{j}^{i}(x),s_{jk}^{i}(x))$ which we shortly write as $%
s(x)=(s_{j}^{i}(x),s_{jk}^{i}(x))=(s_{1}(x),s_{2}(x)),$ using our symbolic
notation in Section 2. We define $g_{1}(x)\overset{def}{=}%
F_{1}(s_{1}(x),s_{2}(x))$ $=s_{1}(x)^{T}s_{1}(x)$ and $g_{2}(x)\overset{def}{%
=}F_{2}(s_{1}(x),s_{2}(x))=s_{1}(x)^{-1}s_{2}(x)$, that is

\begin{equation}
g(x)\overset{def}{=}(g_{1}(x),g_{2}(x))\overset{def}{=}%
F((s_{1}(x),s_{2}(x))=F(s(x))
\end{equation}%
We claim that (20) does not depend on the section $s(x).$ Indeed, if $t(x)$
is another such section, then $t(x)\circ s(x)^{-1}\in \mathcal{G}%
_{2}^{e,e}=\varepsilon O(n)$ since the coordinate system around $e$ is
regular, and (14) implies $F(s(x))=F(t(x))$. Now a coordinate change $%
(U,x^{i})\rightarrow (V,y^{i})$ transforms the components of the section $%
s(x)$ as

\begin{equation}
(s_{j}^{i}(y),s_{jk}^{i}(y))=(s_{j}^{i}(x),s_{jk}^{i}(x))\ast (\frac{%
\partial x^{i}}{\partial y^{j}},\frac{\partial ^{2}x^{i}}{\partial
y^{j}\partial y^{k}})
\end{equation}%
or shortly \ \ $s(y)=s(x)\ast (\frac{\partial x}{\partial y})$ where $\ast $
denotes the group operation of $G_{2}(n)$ defined by (3). From (17), (18)
and (21) we deduce

\begin{eqnarray}
g_{ij}(y) &=&g_{ab}(x)\frac{\partial x^{a}}{\partial y^{i}}\frac{\partial
x^{b}}{\partial y^{j}} \\
g_{jk}^{i}(y) &=&g_{bc}^{a}(x)\frac{\partial x^{b}}{\partial y^{j}}\frac{%
\partial x^{c}}{\partial y^{k}}\frac{\partial y^{i}}{\partial x^{a}}+\frac{%
\partial y^{i}}{\partial x^{a}}\frac{\partial ^{2}x^{a}}{\partial
y^{i}\partial y^{j}}
\end{eqnarray}%
The transformation law (22), (23) satisfies the group law by our remark at
the end of Section 2. So we defined a second order geometric object $g$ on $%
M $ with components $(g_{ij}(x),g_{jk}^{i}(x))$ on $(U,x^{i})$ subject to
the transformation laws (22), (23). Observe that $g$ is constructed from the 
$2$-arrows of $\mathcal{G}_{2}$ using an invariance condition, that is, $g$
has no seperate presence of its own, at least from the present standpoint.
We considered \textit{right }cosets in (14) and fixed the \textit{target }in
the definition of $s$ in order to deal with $g_{ij}$ rather than $g^{ij},$
but such choices are not much relevant for the theory.

At first sight, it seems that $g$ is made up of two seperate geometric
objects, a metric defined by (22) and a torsionfree affine connection
defined by (23). This peculiarity is due to the splitting of $G_{1}(n)$
inside $G_{2}(n).$ However, a closer look reveals that (22) and (23) are
related in a subtle way. To see this, some $f^{p,e}\in \mathcal{G}_{2}^{p,e}$
defines some reqular coordinates $(U,x^{i})$ around $p$ by pulling back the
one at $e$ to $p$. If we choose our section $s$ with $s(p)=f^{p,e}$ and
define $g(x)$ using this particular $s(x)$, we find

\begin{equation}
g_{ij}(p)=\delta _{ij}\text{ \ , \ \ }g_{jk}^{i}(p)=0
\end{equation}

Now (24) shows that (22) and (23) live together, justifying our notation $g,$
and the auxiliary object defined by (23) is not far from the Levi-Civita
connection (we use the terms \textquotedblleft Levi-Civita
connection\textquotedblright\ and \textquotedblleft Christoffel
symbols\textquotedblright\ synonomously). At this point, it is natural to
ask why we work with some auxiliary objects which imitate the Levi-Civita
connection but not work with the Levi-Civita connection itself. Propositions
11, 12 will give a rather unexpected answer to this fair question.

Henceforth in this note a Riemannian structure means a structure as defined
above (we could not find a better name!). Since Definition 4 already
incorporates the metric (22), Riemannian geometry $RG$ \textquotedblleft
includes\textquotedblright\ metric Riemannian geometry $MRG.$ Our purpose in
this note is \textit{not }to show that the inclusion\ $MRG\subset RG$ is
proper (meaning that $RG$ gives new results in $MRG$, we do not know this),
but to show that $RG$ generalizes in a straightforward way to
pre-homogeneous geometries in such a way that one can completely avoid
torsion, covariant differentiation and the Levi-Civita connection. This
generalization will be based on Lie's theorems as we will see in Section 4.

Now let $(U,x^{i}),$ $(V,y^{i})$ be two coordinate patches on $M.$ Using
(22), (23), it is now easy to show that some $2$-arrow $(x^{i},y^{i},\phi
_{j}^{i}(x,y),\phi _{jk}^{i}(x,y))$ of $\mathcal{U}_{2}$ with source in $U$
and target in $V$ belongs to $\mathcal{G}_{2}$ if and only if it satisfies

\begin{eqnarray}
g_{ij}(x) &=&g_{ab}(y)\phi _{i}^{a}(x,y)\phi _{j}^{b}(x,y) \\
g_{jk}^{a}(x)\phi _{a}^{i}(x,y) &=&g_{bc}^{i}(y)\phi _{j}^{b}(x,y)\phi
_{k}^{c}(x,y)+\phi _{jk}^{i}(x,y)  \notag
\end{eqnarray}

(25) gives a set of equations which define the submanifold $\mathcal{G}_{2}$ 
$\subset \mathcal{U}_{2}$ locally. In short, $\mathcal{G}_{2}$ consists of
all $2$-arrows which preserve the geometric object $g.$ Note that if $(\phi
_{j}^{i}(x,y),\phi _{jk}^{i}(x,y))$ and $(\phi _{j}^{i}(x,y),\overline{\phi }%
_{jk}^{i}(x,y))$ both solve (25), then $\phi _{jk}^{i}(x,y)=\overline{\phi }%
_{jk}^{i}(x,y).$

Even though $\mathcal{G}_{2}$ looks like a purely geometric object at first
sight, it is actually a nonlinear system of $PDE$'s made up of inital
conditions which are its $2$-arrows. More precisely, let $f^{p,q}\in 
\mathcal{G}_{2},$ choose coordinates $(U,x^{i}),$ $(V,y^{i})$ around $p,q$
and write $f^{p,q}=(\overline{x}^{i},\overline{y}^{i},\phi _{j}^{i}(%
\overline{x},\overline{y}),\phi _{jk}^{i}(\overline{x},\overline{y})).$
Clearly, the substitution of the components of $(\overline{x}^{i},\overline{y%
}^{i},\phi _{j}^{i}(\overline{x},\overline{y}),\phi _{jk}^{i}(\overline{x},%
\overline{y}))$ into (25) gives an identity since $f^{p,q}\in \mathcal{G}%
_{2}.$ Suppose there exists a local diffeomorphism $f:U\rightarrow V$ which
satisfies the inital condition $f^{i}(\overline{x})=\overline{y}^{i},\frac{%
\partial f^{i}}{\partial x^{j}}(\overline{x})=\phi _{j}^{i}(\overline{x},%
\overline{y}),$ $\frac{\partial ^{2}f^{i}}{\partial x^{j}\partial x^{k}}(%
\overline{x})=\phi _{jk}^{i}(\overline{x},\overline{y})$, and the
substitution $f^{i}(x)=y^{i},$ $\frac{\partial f^{i}}{\partial x^{j}}%
(x)=\phi _{j}^{i}(x,y)$, $\frac{\partial ^{2}f^{i}}{\partial x^{j}\partial
x^{k}}(x)=\phi _{jk}^{i}(x,y)$ satisfies (25) identically for all $x$, that
is, all $2$-arrows of $f$ belong to $\mathcal{G}_{2}$. In this case, we call 
$f$ a local solution of $\mathcal{G}_{2}$ with the inital condition $f^{p,q}$%
. Clearly, a local solution $f$ satisfies all the initial conditions defined
by its $2$-arrows. A global solution of $\mathcal{G}_{2}$ is a
diffeomorphism $f\in Diff(M)$ such that $f_{\mid U}$ is a local solution for
all coordinate patches $(U,x^{i})$ on $M.$ Now $\mathcal{G}_{2}$ admits one
global solution, namely $id,$ because $j_{2}(id)^{p,p}\in $ $\mathcal{G}%
_{2}^{p,p}$ for all $p$. However, $\mathcal{G}_{2}$ may not admit any other
local solutions. The other extreme is a fundamental concept.

\begin{definition}
$\mathcal{G}_{2}$ is completely integrable if

$i)$ All $2$-arrows of $\mathcal{G}_{2}$ integrate to local solutions

$ii)$ A local solution is uniquely determined on its domain by any of its $2$%
-arrows.
\end{definition}

Observe that $2$-arrows of a local solution are determined by its $1$-arrows
in view of the splitting $\varepsilon $ in (19). Thus a local solution is
determined also by any of its $1$-arrows. We will see in Section 5 that $%
i)\Rightarrow ii)$. The reason is that $\mathcal{G}_{2}$ has the property $%
\mathcal{G}_{2}\simeq \mathcal{G}_{1}=\pi \mathcal{G}_{2}.$ The complete
integrability of $\mathcal{G}_{2}$ is a local condition which can be checked
on coordinate patches $(U,x^{i}),$ because if all \textquotedblleft
short\textquotedblright\ $2$-arrows of $\mathcal{G}_{2}$ integrate to local
solutions, then \textit{all }$2$-arrows of $\mathcal{G}_{2}$ integrate to
local solutions. This fact is easily shown as in the proof of Proposition
7.5 in [2]. If $\mathcal{G}_{2}$ is completely integrable, then its local
solutions form a pseudogroup on $M$ because $2$-arrows of $\mathcal{G}_{2}$
are closed under composition and inversion by the definition of a groupoid.
We will denote this pseudogroup by $G$ and its restriction to some $%
(U,x^{i}) $ by $G_{\mid U}.$

Now we want to linearize the $PDE$ $\mathcal{G}_{2}$, which amounts to
defining its algebroid $\mathfrak{G}_{2}.$ First, we recall that the
algebroid of $\mathcal{U}_{2}$ is the vector bundle $J_{2}T\rightarrow M$
whose fiber over $p\in M$ consists of $2$-jets of vector fields at $p.$ So a
section $X_{2}$ of $J_{2}T\rightarrow M$ (with an abuse of notation, we will
write $X_{2}\in J_{2}T)$ is of the form $%
(X^{i}(x),X_{j}^{i}(x),X_{jk}^{i}(x))$ over $(U,x^{i}).$ There is a bracket $%
[$ $,$ $]$ defined on the sections of $J_{2}T\rightarrow M,$ called the
Spencer bracket, which turns $J_{2}T\rightarrow M$ into an algebroid. To
define $[$ $,$ $],$ recall the Spencer operator $D:J_{3}T\rightarrow
J_{2}T\otimes T^{\ast }$ locally given by $%
(X^{i}(x),X_{j}^{i}(x),X_{jk}^{i}(x),X_{mjk}^{i}(x))\rightarrow (\partial
_{j}X^{i}(x)-X_{j}^{i}(x),\partial _{j}X_{k}^{i}(x)-X_{jk}^{i}(x),\partial
_{m}X_{jk}^{i}(x)-X_{mjk}^{i}(x)).$ We have the algebraic bracket $\{$ , $%
\}_{p}:(J_{3}T)_{p}\times (J_{3}T)_{p}\rightarrow (J_{2}T)_{p}$ whose local
formula is obtained by differentiating the usual formula for the bracket of
two vector fields three times at $p$ and replacing derivatives with jet
variables. Clearly, $\{$ $,$ $\}_{p}$ extends to sections of $%
J_{3}T\rightarrow M$ which we denote by $\{$ $,$ $\}$. Now if $%
X_{2},Y_{2}\in J_{2}T$, their Spencer bracket is defined by

\begin{equation}
\lbrack X_{2},Y_{2}]\overset{def}{=}\{X_{3},Y_{3}%
\}+i_{X_{0}}D(Y_{3})-i_{Y_{0}}D(X_{3})
\end{equation}

In (26), $X_{3},Y_{3}$ are arbitrary lifts of $X_{2},Y_{2}$ to $J_{3}T,$ $%
X_{0}=\pi X_{2},$ $Y_{0}=\pi Y_{2}$ where $\pi :J_{2}T\rightarrow J_{0}T=T$
is the projection and $i_{Z}$ denotes contraction with respect to the vector
field $Z\in J_{0}T.$ The bracket $[X_{2},Y_{2}]$ is independent of the lifts 
$X_{3},Y_{3}.$ If $X_{2}=(X^{i}(x),X_{j}^{i}(x),X_{jk}^{i}(x))$ and $%
Y_{2}=(Y^{i}(x),Y_{j}^{i}(x),Y_{jk}^{i}(x))$, we compute

\begin{eqnarray}
\lbrack X_{2},Y_{2}]^{i} &=&X^{a}\partial _{a}Y^{i}-Y^{a}\partial _{a}X^{i}
\\
\lbrack X_{2},Y_{2}]_{j}^{i}
&=&X_{j}^{a}Y_{a}^{i}-Y_{j}^{a}X_{a}^{i}+X^{a}\partial
_{a}Y_{j}^{i}-Y^{a}\partial _{a}X_{j}^{i}  \notag \\
\lbrack X_{2},Y_{2}]_{jk}^{i}
&=&X_{jk}^{a}Y_{a}^{i}+X_{j}^{a}Y_{ka}^{i}+X_{k}^{a}Y_{aj}^{i}-Y_{jk}^{a}X_{a}^{i}-Y_{j}^{a}X_{ka}^{i}-Y_{k}^{a}X_{aj}^{i}
\notag \\
&&+X^{a}\partial _{a}Y_{jk}^{i}-Y^{a}\partial _{a}X_{jk}^{i}  \notag
\end{eqnarray}

Now (27) shows that the projection maps $J_{2}T\rightarrow J_{1}T\rightarrow
J_{0}T$ preserve brackets. Sometimes we will use the same notation $[$ $,$ $%
] $ for all these brackets. Of course, $[$ $,$ $]$ has all the properties
one expects from a bracket (see [28], [25] for further details). Also, we
have the prolongation map $j_{2}:J_{0}T\rightarrow J_{2}T$ defined locally
by $X^{i}(x)\rightarrow (X^{i}(x),\frac{\partial X^{i}}{\partial x^{j}}(x),%
\frac{\partial ^{2}X^{i}}{\partial x^{k}\partial x^{j}}(x))$ and (27) shows
that $[ $ $,$ $]$ respects prolongation, that is, $%
j_{2}[X,Y]=[j_{2}X,j_{2}Y].$

Now, rather than defining $\mathfrak{G}_{2}$ abstractly, we will take a
shortcut following [25] and derive the defining equations of $\mathfrak{G}%
_{2}$ from (25). This method allows one to do explicit computations in
coordinates, but leaves the fundamental relation between $\mathfrak{G}_{2},$ 
$\mathcal{G}_{2}$ which we will need in Section 4 in dark, as we will see.
So we substitute $y^{i}=x^{i}+tX^{i}(x)$ into $%
g(y)=(g_{ij}(y),g_{jk}^{i}(y)) $ in (25), substitute $\phi
_{j}^{i}(x,y)=\delta _{j}^{i}+tX_{j}^{i}(x),$ $\phi
_{jk}^{i}(x,y)=tX_{jk}^{i}(x)$ in (25) and differentiate the resulting
equations with respect to $t$ at $t=0.$ The result is

\begin{eqnarray}
0 &=&X^{a}(x)\partial
_{a}g_{ij}(x)+g_{ai}(x)X_{j}^{a}(x)+g_{aj}(x)X_{i}^{a}(x) \\
0
&=&X_{jk}^{i}(x)+g_{ak}^{i}(x)X_{j}^{a}(x)+g_{aj}^{i}(x)X_{k}^{a}(x)-g_{jk}^{a}(x)X_{a}^{i}(x)+X^{a}(x)\partial _{a}g_{jk}^{i}(x)
\notag
\end{eqnarray}

Now (28) defines a bundle of vectors $\mathfrak{G}_{2}\rightarrow M$ whose
fiber over $p\in (U,x^{i})$ consists of those points $%
(X^{i}(p),X_{j}^{i}(p),X_{jk}^{i}(p))$ of $J_{2}T$ which satisy (28).
Henceforth, we will omit the variable $x$ in (28) and use the same notation
for points $X_{2}\in J_{2}T$ and sections $X_{2}\in J_{2}T$ as before. Since 
$\dim \mathcal{G}_{2}=2\dim M+\dim O(n)$, we have $\dim (\mathfrak{G}%
_{2})_{p}=\dim M+\dim O(n)=n+\frac{n(n-1)}{2}=\frac{n(n+1)}{2}$ for all $%
p\in M$ where $(\mathfrak{G}_{2})_{p}$ denotes the fiber of $\mathfrak{G}%
_{2}\rightarrow M$ over $p.$ So $\mathfrak{G}_{2}\rightarrow M$ is a vector
bundle of rank $\frac{n(n+1)}{2}$. This fact can be checked also directly
from (28) using regular coordinates (see below). The fundamental fact is
that the sections of $\mathfrak{G}_{2}\rightarrow M$ are closed with respect
to the Spencer bracket $[$ $,$ $].$ This follows from the theory ([28],
[25]), but can be checked directly using (27) and (28). Thus we obtain the
algebroid $\mathfrak{G}_{2}\rightarrow M$ and the infinitesimal version of
(25):

\begin{equation}
\begin{array}{ccc}
J_{2}T & \overset{\pi }{\longrightarrow } & J_{1}T \\ 
\cup &  & \cup \\ 
\mathfrak{G}_{2} & \overset{\pi }{\longrightarrow } & \mathfrak{G}_{1}%
\end{array}%
\end{equation}%
where $\mathfrak{G}_{2}\simeq \pi \mathfrak{G}_{2}\overset{def}{=}\mathfrak{G%
}_{1}.$ The splitting $\varepsilon \overset{def}{=}(\pi _{\mid \mathfrak{G}%
_{2}})^{-1}$ amounts to expressing $X_{jk}^{i}(x)$ in (28) in terms of the
lower order terms. Since $[\mathfrak{G}_{2},\mathfrak{G}_{2}]\subset 
\mathfrak{G}_{2}$, we have $[\varepsilon X_{1},\varepsilon
Y_{1}]=\varepsilon \lbrack X_{1},Y_{1}],$ $X_{1},X_{2}\in \mathfrak{G}_{1}$.
Note that the Lie algebra of sections of $\mathfrak{G}_{2}\rightarrow M$\
(which we denoted by $\mathfrak{G}_{2})$ is an \textit{infinite }dimensional
Lie algebra.

By construction, $\mathfrak{G}_{2}\rightarrow M$ is a linear system of $PDE$%
's. A local solution is a vector field $X=X^{i}(x)$ such that the
substitutions $\frac{\partial X^{i}}{\partial x^{j}}=X_{j}^{i},$ $\frac{%
\partial ^{2}X^{i}}{\partial x^{j}\partial x^{k}}=X_{jk}^{i}$ identically
satisfy (28), that is, the prolongation $j_{2}(X)$ belongs to $\mathfrak{G}%
_{2}$. So the fibers of $\mathfrak{G}_{2}\rightarrow M$ consist of initial
conditions. The zero section is of course a local solution, but $\mathfrak{G}%
_{2}\rightarrow M$ may not admit any other local solutions. On the other
extreme, we have

\begin{definition}
$\mathfrak{G}_{2}\rightarrow M$ \ is completely integrable if

$i)$ For any inital condition $\xi \in (\mathfrak{G}_{2})_{p},$ there exists
a local solution around $p$ satisfying this inital condition.

$ii)$ A local solution is uniquely determined on its domain by any of its $2$%
-jets.
\end{definition}

Again, a local solution is actually determined by any of its $1$-jets in
view of the splitting $\varepsilon $ in (29) and $i)\Rightarrow ii)$ for the
same reason: $\mathfrak{G}_{2}\simeq \mathfrak{G}_{1}=\pi \mathfrak{G}_{2}$
as we will see in Section 5$.$ Now the fundamental fact is that the local
solutions around $p$ satisfying the initial conditions $(\mathfrak{G}%
_{2})_{p}$ are closed with respect to the bracket $\ [$ $,$ $]_{J_{0}T}$,
because $[\mathfrak{G}_{2},\mathfrak{G}_{2}]_{J_{2}T}\subset \mathfrak{G}%
_{2} $ and $[$ $,$ $]$ respects prolongation. However, we may have to
restrict the domains of the local solutions to compute their bracket. Taking
the germs of solutions at $p$ as the stalk of a sheaf at $p$, we obtain a
coherent sheaf of Lie algebras defined on $M$, but we will prefer to work
with the more intuitive presheaf in the next section. Clearly all stalks of
this sheaf are isomorphic (see below). If $\mathfrak{G}_{2}$ is completely
integrable, what possible choices do we have for this Lie algebra?

It is a fundamental fact that $\mathcal{G}_{2}$ is completely integrable $%
\Leftrightarrow $ $\mathfrak{G}_{2}$ is completely integrable. In fact, $%
\Rightarrow $ is Lie's first theorem, $\Leftarrow $ is Lie's second theorem
and $\Leftrightarrow $ amounts to constructing the exponential map $\exp :%
\mathfrak{G}_{2\mid U}$ $\rightarrow $ $\mathcal{G}_{2\mid U}$ as envisioned
by Lie and will be sketched in the next section.

\section{Lie's theorems, completeness and uniformization}

Our purpose in this section is to outline a Lie theoretic derivation of the
well known uniformization theorem of space forms in $MRG$. We will try to
avoid the fundamental concepts of $MRG$ as much as possible in order to
emphasize that the present approach generalizes in a quite straightforward
way to all prehomogeneous structures where the main tools of $MRG$ will not
be readily available.

Suppose that $\mathfrak{G}_{2}$ is completely integrable. Let $\mathfrak{X}$
denote the Lie algebra of vector fields on $M$ and $\mathfrak{X}_{g}\subset 
\mathfrak{X}$ the vector space of global solutions of $\mathfrak{G}_{2}.$ We
have seen that $\mathfrak{X}_{g}$ is a Lie algebra, but we may have $\dim 
\mathfrak{X}_{g}=0.$ Let $\mathfrak{X}_{g}(U)$ denote the local solutions of 
$\mathfrak{G}_{2}$ on $(U,x^{i}).$ We fix some $p\in U$ and define

\begin{eqnarray}
j_{2}(p) &:&\mathfrak{X}_{g}(U)\longrightarrow (\mathfrak{G}_{2})_{p} \\
&:&X=X^{i}(x)\longrightarrow (j_{2}X)(p)=(X^{i}(p),\frac{\partial X^{i}}{%
\partial x^{j}}(p),\frac{\partial ^{2}X^{i}}{\partial x^{k}\partial x^{j}}%
(p))  \notag
\end{eqnarray}

As we remarked above, $j_{2}X=\varepsilon (j_{1}X),$ $X\in \mathfrak{G}_{2}.$
Recalling the algebraic bracket $\{$ $,$ $\}_{p}:$ $(J_{3}T)_{p}\times
(J_{3}T)_{p}\rightarrow (J_{2}T)_{p},$ we have

\begin{equation}
\left( j_{2}[X,Y]\right) (p)=\{\varepsilon ((j_{1}(X)(p)),\varepsilon
((j_{1}(Y)(p))\}_{p}\text{ \ \ }X,Y\in \mathfrak{X}_{g}
\end{equation}%
because the last two terms on the $RHS$ of (26) vanish on solutions by the
definition of the Spencer operator and the algebraic bracket coincides with
the Spencer bracket. Note that the $LHS$ of (31) needs complete
integrability for its definition whereas its $RHS$ is still defined if we
replace $j_{1}(X)(p),$ $j_{1}(Y)(p)$ with arbitrary $\xi _{1},\eta _{1}\in (%
\mathfrak{G}_{2})_{p}$. Since $[$ $,$ $]$ respects prolongation, $j_{2}(p)$
is a homomorphism of Lie algebras. If we choose $U$ also simply connected,
then any local solution with some initial condition in $(\mathfrak{G}%
_{2})_{p}$ extends uniquely to $U$ (see below) and $j_{2}(p)$ becomes an
isomorphism of Lie algebras. Therefore $\dim \mathfrak{X}_{g}(U)=\frac{n(n+1)%
}{2}$ if $U$ is simply connected which we will assume below.

Let $\mathfrak{X}_{g}(U,p)\subset $ $\mathfrak{X}_{g}(U)$ denote the
solutions which vanish at $p.$ Now $\mathfrak{X}_{g}(U,p)$ is a subalgebra
of dimension $\frac{n(n-1)}{2}=\dim O(n).$ We call $\mathfrak{X}_{g}(U,p)$
the stabilizer subalgebra at $p.$ In coordinates, the definition of $%
\mathfrak{X}_{g}(U,p)$ amounts to setting $X^{i}(p)=0$ in (28). If we choose
our coordinates regular at $p,$ the first formula in (28) shows that $%
X_{j}^{i}(p)$ is skewsymmetric and the second formula in (28) gives $%
X_{kj}^{i}(p)=0$ so that $\mathfrak{X}_{g}(U,p)$ can be identified (not
canonically!) with $\mathfrak{o}(n).$ The isomorphism (30) identifies $%
\mathfrak{X}_{g}(U,p)$ with its image which is, of course, the splitting in
(6) on the level of Lie algebras. Thus $\mathfrak{X}_{g}(U,p)\simeq \mathcal{%
L}(\mathcal{G}_{2}^{p,p})$ $=$ the Lie algebra of the vertex group $\mathcal{%
G}_{2}^{p,p}.$ Note again that the definition of $\mathfrak{X}_{g}(U,p)$
needs complete integrability whereas the definition of $\mathcal{L}(\mathcal{%
G}_{2}^{p,p})$ does not. The bracket of $\mathfrak{X}_{g}(U,p)$ can be seen
from (27) and reduces to the bracket of two orthogonal matrices in regular
coordinates. The key fact is that the components $X^{i}$ pair with
differential terms in (27), (28) so that the substitution $X^{i}=0$ turns
everything to algebra \textquotedblleft in the vertical
direction\textquotedblright . Clearly the choice of $U$ is arbitrary in
these arguments and we have the same local scenario everywhere on $M.$

We will now sketch the construction of the exponential map. We assume that $%
\mathfrak{G}_{2}$ is completely integrable. We define $i\overset{def}{=}%
j_{2}(p)^{-1}$ and fix some $2$-arrow $f^{p,q}\in \mathcal{G}_{2}$ where $%
q\in U$ is arbitrary. For $X_{2}\in (\mathfrak{G}_{2})_{p},$ let $%
f(iX_{2})(t,x)$ denote the $1$-parameter group of local diffeomorphism
generated by $iX_{2}$ such that $f(iX_{2})(0,p)=p.$ We consider the equation 
$f(iX_{2})(t,p)=q$ in the unknowns $t$ and $X_{2}.$ Since $\mathfrak{G}%
_{0}=J_{0}T=T,$ for any two solutions $t,X_{2}=(X^{i},X_{j}^{i},X_{kj}^{i})$
and $\overline{t},\overline{X}_{2}=(\overline{X}^{i},\overline{X}_{j}^{i},%
\overline{X}_{kj}^{i})$ we have $t=\overline{t}$ and $X^{i}=\overline{X}%
^{i}. $ With these unknowns solved uniquely, we now have the freedom for $%
X_{j}^{i} $ to account for the $1$-arrow $\pi f^{p,q}\in \mathcal{G}_{1}.$
We should make this choice such that $X_{2}=(X^{i},X_{j}^{i},X_{kj}^{i}),$
now determined by $X_{1}=(X^{i},X_{j}^{i}),$ will give the desired $2$-arrow 
$f^{p,q}$. It is a remarkable and nontrivial fact that this can be done only
in one way. In fact, it turns out that the equation $\left[
j_{2}f(iX_{2})(t,x)\right] ^{p,q}=f^{p,q}$ uniquely determines $t$ and $%
X_{2} $ in such a way that all $2$-arrows $\left[ j_{2}f(iX_{2})(t,x)\right]
^{r,s} $ with $r$ close to $p$ belong to $\mathcal{G}_{2},$ that is, $%
j_{2}f(iX_{2})(t,x)$ is a local solution of $\mathcal{G}_{2}$ satisfying the
initial condition $f^{p,q}$. Therefore, all $2$-arrows of $\mathcal{G}_{2}$
starting from $p$ (or ending at $p)$ integrate to local solutions. It
follows that all $2$-arrows of $\mathcal{G}_{2}$ inside $U$ also integrate
to local solutions because any $2$-arrow inside $U$ is a composition (not
uniquely) of two such $2$-arrows. So we conclude that $\mathcal{G}_{2}$ is
completely integrable. The converse follows along the same lines and amounts
to the trick of deriving (28) from (25).

So we have

\begin{proposition}
(Lie's 1st and 2nd theorems) $\mathcal{G}_{2}$ is completely integrable if
and only if $\mathfrak{G}_{2}$ is completely integrable
\end{proposition}

Now we want to globalize Lie's theorems.

If $\gamma $ is any continuous path from $a$ to $b$ with $a\in U$, then we
can continue $\mathfrak{X}_{g}(U)$ uniquely along $\gamma $ like analytic
continuation. However, we may not be able to continue indefinitely as the
local solutions may not approach a definite limit as we approach some point
on the path.

\begin{definition}
\ Suppose $\mathfrak{G}_{2}$ is completely integrable. Then $\mathfrak{G}%
_{2} $ is complete if all local solutions can be continued (necessarily
uniquely) indefinitely along all paths.
\end{definition}

Observe the weakness of Definition 8: it needs complete integrability to
define completeness whereas geodesic completeness of the metric needs no
assumptions on the curvature. If $M$ is compact, then $\mathfrak{G}_{2}$ is
complete. The proof is identical to the proof of Lemma 7.3 in [2]. If $%
\mathfrak{G}_{2}$ is complete and $M$ is simply connected, then the standard
monodromy argument shows that any local solution globalizes uniquely to a
global solution. Even if $M$ is not simply connected, local solutions may
globalize. In this case we call $\mathfrak{G}_{2}$ globalizable. For local
Lie groups globalizability is defined and studied first in [19]. This is a
subtle concept which we will not touch here (see [2] for a cohomological
obstruction to globalizability for $m=0$). However, the algebroid $\mathfrak{%
G}_{2}$ always lifts to an algebroid $\widetilde{\mathfrak{G}_{2}}$ on $%
\widetilde{M}$ where $\rho :\widetilde{M}\rightarrow M$ is the universal
covering space. If $\mathfrak{G}_{2}$ is complete, then so is $\widetilde{%
\mathfrak{G}_{2}}$ (see the proof of Proposition 7.4 in [2]). Since $%
\widetilde{M}$ is simply connected, $\widetilde{\mathfrak{G}_{2}}$
globalizes on $\widetilde{M}.$ We will omit the rather straightforward
details of these arguments. Clearly $\widetilde{\mathfrak{G}_{2}}$ is
locally \textquotedblleft the same\textquotedblright\ as $\mathfrak{G}_{2}.$
We denote the Lie algebra of global solutions of $\widetilde{\mathfrak{G}_{2}%
}$ by $\widetilde{\mathfrak{X}_{g}}.$ By the construction of $\widetilde{%
\mathfrak{X}_{g}}$, any $X\in \widetilde{\mathfrak{X}_{g}}$ is globally
determined by its $2$-jet (or $1$-jet) at any point $p\in \widetilde{M}$ and
therefore $\dim \widetilde{\mathfrak{X}_{g}}=\frac{n(n+1)}{2}.$ The
subalgebra $\widetilde{\mathfrak{X}_{g}}(p)\subset \widetilde{\mathfrak{X}%
_{g}}$ consisting of the vector fields vanishing at $p$ is isomorphic to $%
\mathcal{L}(\mathcal{G}_{2}^{p,p})\simeq \mathfrak{o}(n).$ Indeed, any
statement about $\mathfrak{X}_{g}(U)$ has a global analog for $\widetilde{%
\mathfrak{X}_{g}}.$

We will now repeat the above construction by replacing $\mathfrak{G}_{2}$
with $\mathcal{G}_{2}$ assuming complete integrability. There are no new
ideas involved and we will be very brief. First, by choosing $U$ simply
connected, we may assume that the local solutions of $\mathcal{G}_{2}$
inside $U$, that is, the elements of the pseudogroup $G_{\mid U}$ are all
defined on $U.$ We define the stabilizer $H_{p}\overset{def}{=}\{f\in
G_{\mid U}\mid f(p)=p\}.$ The exponential map integrates the Lie algebra $%
\mathfrak{X}_{g}(U)$ to $G_{\mid U}$ and the stabilizer subalgebra $%
\mathfrak{X}_{g}(U,p)$ to $H_{p}\simeq $ $\mathcal{G}_{2}^{p,p}$.

Now a local solution has a unique continuation along a path by translating
the source of its $2$-arrow (or $1$-arrow) along the path as in [2] . We
call $\mathcal{G}_{2}$ complete if indefinite continuations are possible.
Now the key fact is that $\mathcal{G}_{2}$ is complete if and only if $%
\mathfrak{G}_{2}$ is complete. To see this, we first observe that the
\textquotedblleft short\textquotedblright\ $2$-arrows continue indefinitely
if and only if all $2$-arrows continue indefinitely. This is shown easily by
expressing an arbitrary $2$-arrow as a composition of short $2$-arrows. Now
the statement follows from the exponential map. If $M$ is simply connected,
then the pseudogroup $G$ globalizes to a transformation group $G$ on $M$
which acts transitively and effectively on $M.$ Any transformation $f\in G$
is globally determined by any of its $2$-arrows (or $1$-arrows) so the
geometric order of $(G,M)$ is one. Clearly $H_{p}\overset{def}{=}\{f\in
G\mid f(p)=p\}$ is isomorphic to $O(n).$ Observe that the Klein geometry $%
(G,M)$ determines the vertex connection $\{G,M,H\}=\mathbf{R}$ that we
started with by construction. If globalization is not possible on $M$, we
lift the pseudogroup $G$ to a pseudogroup $\widetilde{G}$ on $\widetilde{M}$
which globalizes to a global transformation group $\widetilde{G}$ which acts
affectively and transitively on $\widetilde{M}.$ The local properties of the
Klein geometry $(\widetilde{G},\widetilde{M})$ are \textquotedblleft the
same\textquotedblright\ as the pseudogroup $G_{\mid U}$ (which may be viewed
as a local Lie group in the classical sense if consider $2$-arrows eminating
from some based point as in the construction of the exponential map). Now
another key fact is that the deck transformations $Deck(\widetilde{M})\simeq
\pi _{1}(M)$ belong to $\widetilde{G}$ since they commute with the
projection $\widetilde{M}\rightarrow M$ and they act as a discontinuous
group on $\widetilde{M}.$

Before we state the next proposition, we need one abstract construction
which sheds further light on the above scenario. The covering map $\rho :%
\widetilde{M}\rightarrow M$ pulls back the groupoid $\mathcal{G}_{2}$ to a
groupoid $\rho ^{-1}\mathcal{G}_{2}$ by pulling back its arrows. This
construction does not need complete integrability and works with all Lie
groupoids. In particular, the algebroid $\mathfrak{G}_{2}$ pulls back to $%
\rho ^{-1}\mathfrak{G}_{2}$ which is the algebroid of $\rho ^{-1}\mathcal{G}%
_{2}.$ If $\mathcal{G}_{2}$ $(\mathfrak{G}_{2})$ is completely integrable,
then $\rho ^{-1}\mathcal{G}_{2}$ $(\rho ^{-1}\mathfrak{G}_{2})$ is also
completely integrable. As we remarked above, completeness of $\mathfrak{G}%
_{2}$ implies completeness of $\rho ^{-1}\mathfrak{G}_{2}$. Conversely, $%
\rho ^{-1}\mathfrak{G}_{2}$ is complete if and only if $\rho ^{-1}\mathfrak{G%
}_{2}$ is globalizable (and therefore $\mathfrak{G}_{2}$ is also complete,
compare to [14], Theorem 4.6, pg. 176). Similar statements hold for $%
\mathcal{G}_{2}.$ So we showed above that $\rho ^{-1}\mathfrak{G}_{2}$
globalizes to $\widetilde{\mathfrak{X}_{g}}$ if and only if $\rho ^{-1}%
\mathcal{G}_{2}$ globalizes to $\widetilde{G}.$ By the definition of the
exponential map, $\widetilde{G}$ has $\widetilde{\mathfrak{X}_{g}}$ as its
infinitesimal generators.

To summarize, we state

\begin{proposition}
(Lie's 3rd theorem) $\rho ^{-1}\mathfrak{G}_{2}$ integrates to a Lie algebra
of global vector fields $\widetilde{\mathfrak{X}_{g}}$ on $\widetilde{M}$ if
and only if $\rho ^{-1}\mathcal{G}_{2}$ integrates to a global
transformation group $\widetilde{G}$ on $\widetilde{M}$. In both cases, $%
\widetilde{G}$ has $\widetilde{\mathfrak{X}_{g}}$ as its infinitesimal
generators. The Klein geometry $(\widetilde{G},\widetilde{M})$ defines the
vertex connection $\mathbf{R}$ and has the above stated properties.
\end{proposition}

Note also that some $f^{p,q}\in \mathcal{G}_{2}^{p,q}$ induces an
isomorphism $f_{\ast }^{p,q}:(\mathfrak{G}_{1})_{p}\rightarrow (\mathfrak{G}%
_{1})_{q}$ (this association does not need complete integrability, see [3])
which lifts to the adjoint action of $\widetilde{G}$ on its Lie algebra $%
\widetilde{\mathfrak{X}_{g}}.$

Now what possibilities do we have for the Klein geometry $(\widetilde{G},%
\widetilde{M})?$ Let $(G_{1},M_{1}),$ $(G_{2},M_{2})$ be two Klein
geometries. A morphism $(G_{1},M_{1})\rightarrow (G_{2},M_{2})$ is a pair $%
(\varphi ,f)$ where $\varphi :$ $G_{1}\rightarrow G_{2}$ is a Lie group
homomorphism and $f:M_{1}\rightarrow M_{2}$ is a smooth map satisfying $%
f(gx)=\varphi (g)f(x).$ We define an isomorphism $(G_{1},M_{1})\simeq
(G_{2},M_{2})$ in the obvious way. Now it is natural to believe that $(%
\widetilde{G},\widetilde{M})$ is isomorphic to one of

\begin{equation}
\text{\ }(G(1),\mathbb{S}^{n})\text{ \ \ \ \ \ \ \ \ \ \ }(G(0),\mathbb{R}%
^{n})\text{ \ \ \ \ \ \ \ \ \ \ (}G(-1),\mathbb{H}^{n})
\end{equation}

We can deduce (32) from the uniformization theorem of $MRG.$ To do this, we
need to show two more facts.

$1)$ $\mathfrak{G}_{2}$ (or $\mathcal{G}_{2})$ is complete $\Leftrightarrow $
the underlying metric is geodesically complete.

$2)$ $\mathfrak{G}_{2}$ (or $\mathcal{G}_{2})$ is completely integrable $%
\Leftrightarrow $ the underlying metric has constant curvature.

We will prove $2)$ in Section 7. However, it is also possible to deduce (32)
from a Lie theoretic statement which therefore implies uniformization
theorem of $MRG.$ To do this, let $\mathcal{A}$ denote the set of all Klein
geometries $(G,G/H)$ satisfying the following properties.

$i)$ $H=O(n)$ and $O(n)\subset G$ is a Lie subgroup.

$ii)$ $G/O(n)$ is simply connected.

$iii)$ $G$ acts effectively on $G/O(n)$ with geometric order $m=1$

$iv)$ $\{G,G/H,H\}=\mathbf{R}$

The requirement $m=1$ in $iii)$ is redundant by $ii)$ and the compactness of 
$O(n).$ With these assumptions, we need to show that any $(G,G/H)\in $ $%
\mathcal{A}$ is isomorphic to one of (32). Note that $(\widetilde{G},%
\widetilde{M})$ satisfies all the requirements.

The above group theoretic statement is actually a statement about Lie
algebra pairs $(\mathfrak{g,o}(n))$ in view of Remark 1 (see also Section
7). We believe that $iv)$ is also redundant and is a consequence.

At any rate, we will state

\begin{proposition}
(Lie's 3rd theorem, refined form) The Klein geometry $(\widetilde{G},%
\widetilde{M})$ in Proposition 9 is isomorphic to one of (32)
\end{proposition}

\QTP{Body Math}
Propositions 9, 10 may seem somewhat surprising at first, but actually they
are quite expected. Indeed, recall the classical formula for the Lie
derivative of a metric: $\mathcal{L}_{X}(g_{jk})=X^{a}\partial
_{a}g_{ij}(x)+g_{ai}(x)\frac{\partial X^{a}}{\partial x^{j}}+g_{aj}(x)\frac{%
\partial X^{a}}{\partial x^{i}}.$ If we replace $(X^{i},\frac{\partial X^{i}%
}{\partial x^{j}})$ by the $1$-jet $X_{1}=(X^{i},X_{j}^{i})\in J_{1}T,$ we
get the $RHS$ of the first formula in (28). Now any linear geometric object
of arbitrary order can be Lie-differentiated with respect to a vector field
as explained in [34]. If we compute the Lie derivative of the second order
geometric object $g$ with respect to a vector field and replace $(X^{i},%
\frac{\partial X^{i}}{\partial x^{j}},\frac{\partial ^{2}X^{i}}{\partial
x^{k}\partial x^{j}})$ by $X_{2}=(X^{i},X_{j}^{i},X_{kj}^{i})=\varepsilon
(X_{1})\in J_{2}T$, we get the formulas on the $RHS$ of (28). So the
solutions of (28) are Killing vector fields for $g$! Since this computation
deduces Killing vector fields from Lie derivative, covariant differentiation
in tensor calculus \textit{must} be a special case of Lie derivative. The
derivation of the formula (21) in [2] shows that this is indeed the case. We
believe that covariant differentiation in tensor calculus owes its existence
to the splitting of (2) in the exceptional case $k=1,$ $n\geq 1.$ We hope
that this crucial point will become more transparent in the next section.

\section{Curvature}

Definitions 5, 6 have a serious deficiency: they are not effective. How do
we decide complete integrability of $\mathfrak{G}_{2},$ $\mathcal{G}_{2}$
from (25), (28)? Our purpose in this section is to define two curvatures $%
\mathfrak{R},$ $\mathcal{R}$ where $\mathfrak{R}=$ $\mathfrak{G}_{2}$%
-curvature and $\mathcal{R}$ $=$ $\mathcal{G}_{2}$-curvature. We will show
that $\mathfrak{G}_{2}$ is completely integrable $\Leftrightarrow \mathfrak{R%
}=0$ and $\mathcal{G}_{2}$ is completely integrable $\Leftrightarrow 
\mathcal{R}=0.$ Therefore the conditions $\mathfrak{R}=0$, $\mathcal{R}=0$
may be taken as the definitions of complete integrability. The main message
here is that curvature is always an obstruction to complete integrability
for pre-homogeneous structures.

We start with (28). We seperate the second formula in (28) from the first
one and rewrite it as an equivalent first order system

\begin{eqnarray}
\partial _{j}X^{i} &=&X_{j}^{i} \\
\partial _{k}X_{j}^{i}
&=&-g_{ak}^{i}X_{j}^{a}-g_{aj}^{i}X_{k}^{a}+g_{jk}^{a}X_{a}^{i}-X^{a}%
\partial _{a}g_{jk}^{i}  \notag
\end{eqnarray}

Now (33) expresses the derivatives of the unknown functions $%
(X^{i},X_{j}^{i})$ in terms of the functions themselves. The find the
integrability conditions of (33), we differentiate the second formula with
respect to $x^{r},$ substitite back from (33) and alternate $r,k.$ The
result is

\begin{equation}
\mathfrak{R}_{kr,j}^{i}\overset{def}{=}\left[ \widehat{\mathfrak{R}}%
_{kr,j}^{i}\right] _{[kr]}=0
\end{equation}%
where $[kr]$ denotes alternation of the indices $k,r$ and

\begin{eqnarray}
&&\widehat{\mathfrak{R}}_{kr,j}^{i}\overset{def}{=}X_{k}^{a}\partial
_{r}g_{aj}^{i}+X^{b}(\partial _{rb}^{2}g_{jk}^{i}+g_{jk}^{a}\partial
_{b}g_{ra}^{i}-g_{ak}^{i}\partial _{b}g_{jr}^{a})  \notag \\
&&+X_{j}^{b}(\partial
_{r}g_{bk}^{i}-g_{ak}^{i}g_{br}^{a})+X_{r}^{b}(\partial
_{b}g_{jk}^{i}-g_{ak}^{i}g_{bj}^{a}+g_{jk}^{a}g_{ba}^{i})  \notag \\
&&-X_{b}^{i}(\partial _{r}g_{jk}^{b}+g_{jk}^{a}g_{ra}^{b})
\end{eqnarray}

\ Observe the occurence of the the Riemann curvature tensor twice in (35) if
we replace $g_{kj}^{i}$ with the Christoffel symbols $\Gamma _{kj}^{i}.$
Using (24), one can show that this fraud Riemann curvature \textit{tensor}
satisfies many identities as the genuine one.

Now, by the well known existence and uniqueness theorem for first order
systems of $PDE$'s with initial conditions (see, for instance, [17],
pg.224-227), if (35) is satisfied identically for all points $%
X_{1}=(X^{i},X_{j}^{i})\in J_{1}T$, then we may choose an arbitrary inital
condition $(X^{i}(p),X_{j}^{i}(p))\in (J_{1}T)_{p}$ and solve (33) uniquely
for a vector field $X^{i}(x)$ defined around $p$ and satisfying $%
(j_{1}X)(p)=(X^{i}(p),X_{j}^{i}(p)).$ Of course, $j_{2}X\in J_{2}T$ but the
problem is that we may not have $j_{2}X\in $ $\mathfrak{G}_{2}$ even if we
choose $(X^{i}(p),X_{j}^{i}(p))\in (\mathfrak{G}_{1})_{p}$! Obviously, there
is a missing integrability condition which should involve differentiation of
the first formula in (28). So we join the first formula of (28) to (33) and
rewrite (33) as

\begin{eqnarray}
\partial _{j}X^{i} &=&X_{j}^{i} \\
0 &=&X^{a}\partial _{a}g_{ij}+g_{ai}X_{j}^{a}+g_{aj}X_{i}^{a}  \notag \\
\partial _{k}X_{j}^{i}
&=&-g_{ak}^{i}X_{j}^{a}-g_{aj}^{i}X_{k}^{a}+g_{jk}^{a}X_{a}^{i}-X^{a}%
\partial _{a}g_{jk}^{i}  \notag
\end{eqnarray}

Clearly, (28) is equivalent to (36). In the old works, (36) is called a
mixed system due to the constraint coming from the second equation. It is
shown in these works that such a system reduces to a first order system
without any constraint after successive prolongations (see, for instance,
[9], [32]). However, such proofs implicitly assume that some rank conditions
are satisfied so that the implicit function theorem is applicable. Around
1970, these rank conditions are organized by D.C.Spencer and his coworkers
into a powerful technique, now called Spencer cohomology. For instance, see
[26], pg 254-255 for a direct proof in coordinates which derives the
constant curvature condition (49) from the formal integrability of (28) and
makes heavy use of Spencer cohomology (this proof assumes $g_{kj}^{i}$ $%
=\Gamma _{kj}^{i}$ and adds the second equation of (28) to the first
equation as a trick, see pg. 251).

We will now derive the missing integrability condition by an elementary
method which avoids Spencer cohomology. (48) below will justify that this
method recovers all the integrability conditions. So we differentiate the
second equation of (36) with respect to $x^{k}$, substitute $\partial
_{k}X_{j}^{a}$ from the third equation, and alternate $k,i$ in $\partial
_{k}X_{i}^{a}.$ Another straightforward computation gives

\begin{equation}
\mathfrak{R}_{ki,j}\overset{def}{=}\left[ \widehat{\mathfrak{R}}_{ki,j}%
\right] _{[ki]}=0
\end{equation}%
where

\begin{eqnarray}
&&\widehat{\mathfrak{R}}_{ki,j}\overset{def}{=}X_{k}^{a}\partial
_{a}g_{ij}+X^{a}\partial _{ka}^{2}g_{ij}+X_{j}^{a}\partial
_{k}g_{ai}^{a}+X_{i}^{a}\partial _{k}g_{aj} \\
&&-X_{j}^{b}g_{ai}g_{bk}^{a}-X_{k}^{b}g_{ai}g_{bj}^{a}+X_{b}^{a}g_{ai}g_{jk}^{b}-X^{b}g_{ai}\partial _{b}g_{jk}^{a}
\notag
\end{eqnarray}

We define the horizontal (over $M$) $2$-form $\mathfrak{R}$ by

\begin{equation}
\mathfrak{R}_{ij}\overset{def}{=}(\mathfrak{R}_{ij,r},\mathfrak{R}%
_{ij,r}^{s})
\end{equation}

We write $\mathfrak{R}=(\mathfrak{R}_{1},\mathfrak{R}_{2})$ and call $%
\mathfrak{R}$ the algebroid curvature. Now a section of the dual bundle $%
(J_{1}T)^{\ast }\rightarrow M$ is locally of the form $\xi ^{1}=(\xi
_{i},\xi _{i}^{j})$ and pairs with a section $(X^{i},X_{j}^{i})$ of $%
J_{1}T\rightarrow M$ linearly to the function $\xi _{a}X^{a}+\xi
_{a}^{b}X_{b}^{a}.$ We denote this pairing by $(X_{1},\xi ^{1}).$ We believe 
$\mathfrak{R}(X_{0},Y_{0})(Z_{1})\in \mathfrak{G}_{1}^{\ast }$ and $%
\mathfrak{R}(X_{0},Y_{0})(\overline{Z}_{1})\in (J_{1}T)^{\ast }$, $%
X_{0},Y_{0}\in J_{0}T=T,$ $Z_{1}\in \mathfrak{G}_{1},$ $\overline{Z}_{1}\in
J_{1}T.$ A direct proof of these statements in coordinates requires
formidable amount of computation. Assuming this for the moment, it follows
that $\left( W_{1},\mathfrak{R}(X_{0},Y_{0})(Z_{1})\right) $ and $\left( 
\overline{W}_{1},\mathfrak{R}(X_{0},Y_{0})(\overline{Z}_{1})\right) $ \ are
functions on $M.$ In particular, the horizontal $2$-form $\left( \overline{Z}%
_{1},\mathfrak{R}(X_{0},Y_{0})(\overline{Z}_{1})\right) $ lives in the
variational complex and descends to a $2$-form on the quotient $J_{1}T/%
\mathfrak{G}_{1}$ if $\mathfrak{G}_{2}$ is completely integrable. At any
rate, $\mathfrak{R}$ is not a tensor but a second order object. We believe
that the component $\mathfrak{R}_{2}$ is seperated from the full curvature $%
\mathfrak{R}$ and tamed into the Riemann curvature tensor in the same way as
(23) is seperated from (22) and tamed into covariant differentiation. We
should recall here that Riemann writes only one formula in his foundational
Habilitationsschrift and the Riemann curvature tensor is introduced later by
others as a part of covariant differentiation and tensor calculus.

We now turn to (25). We first single out the fraud Riemann curvature tensor.

\begin{equation*}
R_{rj,k}^{i}(x)\overset{def}{=}\left[ \partial
_{r}g_{jk}^{i}(x)-g_{rk}^{b}(x)g_{jb}^{i}(x)\right] _{[rj]}
\end{equation*}

We differentiate the second equation in (25) with respect to $x^{r}$
assuming that $y=y(x)$ is a solution, substitute back $\frac{\partial \phi
_{jk}^{i}(x,y)}{\partial x^{r}}$ from this equation and alternate $r,j.$ For
simplicity of notation, we write $\phi _{j}^{i}$ for $\phi _{j}^{i}(x,y)$
and $\overline{\phi }_{j}^{i}$ for $\phi _{j}^{i}(x,y)^{-1}$, keeping in
mind that $\phi _{j}^{i}$ and $\overline{\phi }_{j}^{i}$ depend on both
source and target variables. The result is

\begin{equation}
\mathcal{R}_{rj,k}^{i}\overset{def}{=}R_{ab,c}^{d}(x)\phi _{d}^{i}(\overline{%
\phi })_{r}^{a}(\overline{\phi })_{j}^{b}(\overline{\phi }%
)_{k}^{c}-R_{rj,k}^{i}(y)=0
\end{equation}

The formula (40) is well known from tensor calculus (see [10] for a proof
that (40) implies (49) when $g_{jk}^{i}=\Gamma _{jk}^{i}$). By the same
method above, we now differentiate the first equation of (25) with respect
to $x^{k}$, substitute $\frac{\partial \phi _{i}^{a}(x,y)}{\partial x^{k}}$
from the second equation and alternate $k,j$ in $\frac{\partial \phi
_{j}^{b}(x,y)}{\partial x^{j}}.$ The final result is

\begin{equation}
\mathcal{R}_{kj}^{i}\overset{def}{=}\left[ \widehat{\mathcal{R}}_{kj}^{i}%
\right] _{[kj]}=0
\end{equation}%
where

\begin{eqnarray}
&&\widehat{\mathcal{R}}_{kj}^{i}\overset{def}{=}(\overline{\phi }%
)_{a}^{b}\partial _{k}g_{bj}(x)g^{ai}(y)-\phi _{k}^{a}\phi _{j}^{b}\partial
_{a}g_{cb}(y)g^{ci}(y) \\
&&-(\overline{\phi })_{d}^{e}\phi _{c}^{a}\phi
_{j}^{b}g_{ab}(y)g_{ke}^{c}(x)g^{di}(y)+g_{ab}(y)g_{cd}^{a}(y)\phi
_{k}^{c}\phi _{j}^{b}g^{di}(y)  \notag
\end{eqnarray}

We define

\begin{equation}
\mathcal{R}_{kj}\overset{def}{=}(\mathcal{R}_{kj}^{i},\mathcal{R}_{kj,r}^{i})
\end{equation}%
write $\mathcal{R}=(\mathcal{R}_{1},\mathcal{R}_{2})$ and call $\mathcal{R}$
the groupoid curvature. Observe that $\mathcal{R}$ depends on both source
and target variables $x,y$ and also on $1$-arrows of $\mathcal{G}_{1}.$ For $%
X_{0},Y_{0}\in T_{p},$ we believe that $\mathcal{R}(X_{0},Y_{0},\phi
^{p,q}):(\mathfrak{G}_{1})_{p}\rightarrow (\mathfrak{G}_{1})_{q}$ is a
linear map.

Our purpose is now to derive the Maurer-Cartan $(MC)$ equations. Recall the
base point $e$ in Section 2 which we fixed once and for all with some
regular coordinates around it. The main idea is to \textquotedblleft
factor\textquotedblright\ all $1$-arrows of $\mathcal{G}_{1}$ through the
point $e$ and thus seperate the source and target variables in $\mathcal{R}$
as on pg. 22 of [2]. Now any $\phi ^{p,q}\in \mathcal{G}_{2}^{p,q}$ is of
the form $(\alpha ^{q,e})^{-1}\circ $ $\beta ^{p,e}$ for some $\alpha
^{q,e}, $ $\beta ^{p,e}\in $ $\cup _{x\in M}\mathcal{G}_{2}^{x,e}.$ Clearly, 
$\alpha ^{q,e},$ $\beta ^{p,e}$ are not unique in this factorization because 
$\overline{\beta }^{p,e}=$ $\lambda ^{e,e}\circ \beta ^{p,e}$ and $\overline{%
\alpha }^{p,e}=$ $\lambda ^{e,e}\circ \alpha ^{p,e}$ work too with arbitrary 
$\lambda ^{e,e}\in \mathcal{G}_{2}^{e,e}$. This factorization is the passage
from the principal bundle to the groupoid. Thus the principal bundle $\cup
_{x\in M}\mathcal{G}_{2}^{x,e}$ and the groupoid $\mathcal{G}_{2}$ are
equivalent objects as they determine each other, \textit{as long as we do
not differentiate.} Indeed, in the derivation of $\mathcal{R}$ in (43), we
regard both source and targets of the arrows as variables when we
differentiate, which is not possible in the principal bundle. At this point,
we will leave it to the reader to work out and judge for himself/herself how
the concept of torsion emerges \textit{as a necessity }from this deficiency
of the principal bundle. We hope that the next computation will further
clarify this point.

We need the above factorization only for $\mathcal{G}_{1}.$ We now have $%
\phi _{j}^{i}(x,y)=\alpha _{a}^{i}(y,e)^{-1}$ $\beta _{j}^{a}(x,e)$. For
simplicity, we omit the base point $e$ from our notation and write $\phi
(x,y)=\alpha (x)^{-1}\circ \beta (y)$, that is,

\begin{equation}
\phi (x,y)_{j}^{i}=\alpha ^{-1}(y)_{a}^{i}\beta (x)_{j}^{a}
\end{equation}

We now substitute (44) into (25). The new equations obtained in this way are
of course equivalent to (25). We now repeat the above derivation of $%
\mathcal{R}$ using these new equations. This amounts to substituting (44)
into (35) and (38). Now (40) becomes

\begin{equation}
\left[ \widetilde{\mathcal{R}}_{rj,k}^{i}(x,\beta (x))\right] _{[rj]}-\left[ 
\widetilde{\mathcal{R}}_{rj,k}^{i}(y,\alpha (y))\right] _{[rj]}=0
\end{equation}

where

\begin{eqnarray*}
&&\widetilde{\mathcal{R}}_{rj,k}^{i}(x,\beta (x))\overset{def}{=}%
R_{ab,c}^{d}(x)\beta _{d}^{i}(x)\beta ^{-1}(x)_{r}^{a}\beta
^{-1}(x)_{j}^{b}\beta ^{-1}(x)_{k}^{c} \\
&&\widetilde{\mathcal{R}}_{rj,k}^{i}(y,\alpha (x))\overset{def}{=}%
R_{ab,c}^{d}(y)\alpha _{d}^{i}(y)\alpha ^{-1}(y)_{r}^{a}\alpha
^{-1}(y)_{j}^{b}\alpha ^{-1}(y)_{k}^{c}
\end{eqnarray*}

We now substitute (44) into the peculiar formula (38). A surprising
computation shows that (41) becomes

\begin{equation}
\left[ \widetilde{\mathcal{R}}_{ik,j}(x,\beta (x))\right] _{[ik]}-\left[ 
\widetilde{\mathcal{R}}_{ik,j}(y,\alpha (y))\right] _{[ik]}=0
\end{equation}

where

\begin{eqnarray}
&&\widetilde{\mathcal{R}}_{ik,j}(x,\beta (x))\overset{def}{=}\beta
^{-1}(x)_{i}^{a}\beta ^{-1}(x)_{j}^{b}\beta ^{-1}(x)_{k}^{c}(\partial
_{a}g_{bc}(x)+g_{da}(x)g_{bc}^{d}(x)) \\
&&\widetilde{\mathcal{R}}_{ik,j}(y,\alpha (y))\overset{def}{=}\alpha
^{-1}(y)_{i}^{a}\alpha ^{-1}(y)_{j}^{b}\alpha ^{-1}(y)_{k}^{c}(\partial
_{a}g_{bc}(y)+g_{da}(y)g_{bc}^{d}(y))  \notag
\end{eqnarray}

Therefore, $\mathcal{R}=0$ if and only if

\begin{eqnarray}
\left[ \widetilde{\mathcal{R}}_{rj,k}^{i}(x,\beta (x))\right] _{[rj]}
&=&c_{rj,k}^{i} \\
\left[ \widetilde{\mathcal{R}}_{rj,k}(x,\beta (x))\right] _{[rj]} &=&c_{rj,k}
\notag
\end{eqnarray}%
for some constants $c_{rj,k}^{i},$ $c_{rj,k}.$ Now (48) is the $MC$
equations for the pseodogroup $G$ or the global transformation group $%
\widetilde{G}$ in Section 4 (see [23] for $MC$ equations for more general
pseudogroups than considered here). The recovery of the group shows that our
elementary method gives all the integrability conditions.

Thus we have

\begin{proposition}
The following conditions are equivalent

i) $\mathcal{G}_{2}$ is completely integrable

ii) $\mathcal{R}=0$ on $\mathcal{G}_{1}$

iii) $\mathfrak{G}_{2}$ is completely integrable

iv) $\mathfrak{R}=0$ on $\mathfrak{G}_{1}$

v) $MC$ equations (48) hold
\end{proposition}

Now the first equation of (28) is formally integrable if and only if the
constant curvature condition

\begin{equation}
\mathfrak{R}_{kr,j}^{i}=c(\delta _{k}^{i}g_{jr}-\delta _{r}^{i}g_{jk})
\end{equation}%
holds. This equivalence is well known. As we remarked above, it is proved in
[25], pg. 254-255 and also in [7], Proposition 2.13. Clearly, complete
integrability implies the formal integrability of the first (in fact both)
equation of (28). Now the weaker concept of formal integrability is
equivalent to complete integrability as shown in [13] in the case of
pseudogroups of finite type. Therefore

\begin{proposition}
The conditions of Proposition 11 are equivalent to

i) (28) is formally integrable

ii) (49) holds
\end{proposition}

Propositions 11, 12 have a surprising consequence. Let $\mathcal{G}_{2},$ $%
\widetilde{\mathcal{G}}_{2}$ two Riemannian structures defining the same
metric structure, that is, $\pi \mathcal{G}_{2}=\pi \widetilde{\mathcal{G}}%
_{2}.$ Then $\mathcal{R}=0$ if and only if $\widetilde{\mathcal{R}}=0.$
Indeed, both conditions are equivalent to (49) which can be checked working
with the metric only. Therefore, \textit{as far as (49) is concerned (but
possibly not further!), }all the auxiliary objects (23) are equal and the
Levi-Civita connection has no priority, answering the fair question in
Section 2. We recall here again that Definition 4 excludes nothing from $MRG$
and what it includes is topologically trivial since Riemannian structures
according to Definition 4 are in 1-1 correspondence with reductions of the
structure group $G_{2}(n)$ of the second order principal (co)frame bundle to
the subgroup $\varepsilon O(n)$ which is homotopically equivalent to $O(n).$

We conclude this section by clarifying some ambiguities in [2] which were
not clear to us at the time of writing [2]. If we linearize the groupoid
with defining equations (39) in [2], we arrive at (21) in [2]. The vector
fields which solve (21) should be called right invariant vector fields
according to the convention in [2] but infinitesimal generators according to
this note. No attention is paid to the infinitesimal generators in [2] but
everything is based on left invariant vector fields. On the other hand, it
is the infinitesimal generators which integrate to solutions of (21) whereas
the left invariant vector fields integrate to right local diffeomorphisms
whose $1$-arrows are computed by the formula (45) in [2]. This corresponds
to the well known fact that left invariant vector fields integrate to right
translations and right invariant vector fields integrate to left
translations on a Lie group. Therefore, Definition 4.1, Lemma 4.2 and some
related arguments are not essential for the main purpose of [2] and torsion
can be avoided also in [2] as in this note.

\section{Affine structures}

\begin{definition}
An affine structure on $M$ is a subgroupoid $\mathcal{A}_{2}$ $\subset 
\mathcal{U}_{2}$ with $\{\mathcal{A}_{2}\}=\mathbf{A.}$
\end{definition}

The defining equations of $\mathcal{A}_{2}$ are given by the second formula
in (25) where we should replace $g_{jk}^{i}$ by $\Gamma _{jk}^{i}.$ So $%
\mathcal{A}_{2}$ is nothing but a torsionfree affine connection on the 
\textit{first order }principal (co)frame bundle of $M.$ Clearly, a
Riemannian structure canonically determines an affine structure. Now the
algebroid $\mathfrak{A}_{2}$ is defined by the second formula in (28). All
the constructions and propositions in Sections 3, 4, 5 carry over word by
word. There is only one Klein geometry up to isomorphism in the
uniformization theorem (see [14], [34] for a thorough study) which can be
derived also from a Lie theoretic statement. The $\mathcal{R}_{1}$ and $%
\mathfrak{R}_{1}$ components of $\mathcal{R}$ and $\mathfrak{R}$ disappear.
The $MC$ equations involve only these components. We omit further details as
our main purpose in this note was to use affine structures to emphasize some
facts in Riemannian geometry.

\section{Pre-homogeneous structures}

This note clearly scratches only the tip of an iceberg and leaves many
questions unanswered even in affine geometry let alone Riemannian geometry.
Still, we feel that it may be useful to formulate some open problems
(admittedly by far obvious) for pre-homogeneous structures with the hope
that they may activate some research.

At this stage, it is quite clear how to define a pre-homogeneous structure.
However there is a technical difficulty we should clarify first. Let $G$ be
a connected Lie group, $H\neq \{e\}$ $\subset G$ a discrete subgroup with $%
H\cap Z(G)=\{e\},$ that is, $G$ acts effectively on $M=G/H.$ If $N$ denotes
the base manifold $G$, then $m=0$ for the Klein geometry $(G,N)$ but $m=1$
for $(G,M)$. Thus the geometric order may decrease by one (but not more, see
[3]) when we pass to a covering. The main point here is that geometric order
is a global concept whereas the vertex connection is actually a local
concept. This is clear from the construction of \ the local pseodogroup $%
G_{\mid U}$ in Section 4 which determines the vertex connection. Indeed, we
need not determine all $g\in G$ in the formula $\mathcal{H}_{m}=\{g\}$ in
the Introduction but only those $g$ near the identity. In [3] we defined
also the infinitesimal order $\overline{m}$ of an effective infinitesimal
Klein geometry $(\mathfrak{g},\mathfrak{h})$. As in Remark 1, if $M$ is
simply connected in $(G,M),$ then $m=\overline{m}.$ So the concept needed in
this note is $\overline{m}$ but we based our study on $m$ because it is more
intuititive and geometric than $\overline{m}$ in the same way as the subtle
nonlinear object $\mathcal{G}_{2}$ is easier to grasp geometrically than its
linearization $\mathfrak{G}_{2}.$

Now let $(G,M)$ be a Klein geometry with geometric order $m$ and $M$ simply
connected. As explained in the Introduction, we have the vertex connection $%
\{G,M,H\}$ as a conjugacy class in $G_{m+1}(n).$ Recalling the definition of 
$\{\mathcal{G}_{m+1}\}$ in Section 3, we make

\begin{definition}
A pre-homogeneous structure on $M$ of order $m+1$ is a transitive Lie
subgroupoid $\mathcal{G}_{m+1}\subset \mathcal{U}_{m+1}$ such that $\{%
\mathcal{G}_{m+1}\}=\{G,N,H\}$ for some Klein geometry $(G,N)$ of geometric
order $m.$ The vertex connection $\{G,N,H\}$ is the model for $\mathcal{G}%
_{m+1}.$
\end{definition}

Therefore, Riemannian and affine structures as defined in Sections 3, 4 are
special pre-homogeneous structures. We now check that a parallelizable
manifold as defined in [2] (see Definition 3.1 in [2]) is another special
case. So let $G$ be a \textit{any }Lie transformation group which acts
simply transitively on $M$. With $H_{0}$ as defined in the Introduction, we
have $H_{0}=\{g\}.$ Therefore $m=0$ and all stabilizers are identity. It
follows that any Klein geometry $(G,M)$ determines the same vertex
connection $\{G,M,\{e\}\}$ with the representative $(o,\delta _{j}^{i})$
where $o$ is the origin of $\mathbb{R}^{n}$ and the isomorphism classes of
such simply connected Lie groups is the same as \textit{all }Lie algebras.
Omitting the base point from our notation, the stabilizer $\{id\}$ injects
into $G_{1}(n)$ as $j_{1}(id)=I$ whose conjugacy class is $\{I\}$.
Therefore, for $m=0,$ a pre-homogeneous structure on $M$ is a first order
groupoid $\mathcal{H}_{1}$ on $M$ with vertex groups $\mathcal{H}%
_{1}^{p,p}=j_{1}(id)^{p,p},$ which is equivalent to the Definition 3.1 in
[2]. To construct the geometric object, the right (or left) coset space $%
G_{1}(n)/I=G_{1}(n)$ and we define the components of $F$ by $F_{j}^{i}(a)%
\overset{def}{=}a_{j}^{i},$ $a\in G_{1}(n),$ so that $ab^{-1}=I%
\Leftrightarrow F(a)=F(b)$ as in (14). Note that $F$ is again globally
defined and a polynomal map of degree one. Now we construct the first order
geometric object $\omega $ in the same way as we did $g$ in Section 3: $%
\omega $ has components $\omega _{j}^{i}(x)$ with $\det a\neq 0$ and $\omega
_{j}^{i}(p)=\delta _{j}^{i}$ in regular coordinates. Some $1$-arrow of $%
\mathcal{U}_{1}$ preserves $\omega $ if and only if it belongs $\mathcal{H}%
_{1}$ and we get the defining equations (39) in [2]. Note that we can define 
$\omega _{a}^{i}(x)\overline{\omega }_{j}^{a}(y)\overset{def}{=}\varepsilon
_{j}^{i}(x,y)$ (see pg. 22 of [2]) and start with $\varepsilon _{j}^{i}(x,y)$
as we did in [2].

Natural Definition 14 may seem, it immediately gives rise to some quite
nontrivial questions. Note that we use Klein geometries to define splittings
inside (2) but if can discover some splitting in some ad hoc way, we can
define the groupoid and get started. This is illustrated by the following
derivation which is possibly new. An element of $G_{3}(1)$ is of the form $%
(a_{1},a_{2},a_{3})$ where $a_{i}\in \mathbb{F=R}$ or $\mathbb{C}.$ The
chain rule gives the group operation of $G_{3}(1):$

\begin{eqnarray}
(a_{1},a_{2},a_{3})(b_{1},b_{2},b_{3}) &=&(c_{1},c_{2},c_{3}) \\
c_{1} &=&a_{1}b_{1}  \notag \\
c_{2} &=&a_{1}b_{2}+a_{2}(b_{1})^{2}  \notag \\
c_{3} &=&a_{1}b_{3}\text{ }+3a_{2}b_{1}b_{2}+a_{3}(b_{1})^{3}  \notag
\end{eqnarray}

We define the map $\varepsilon :G_{2}(1)\rightarrow G_{3}(1)$ by $%
\varepsilon (a_{1},a_{2})=(a_{1},a_{2},\frac{3}{2}\frac{(a_{2})^{2}}{a_{1}})$
and using (50) we check that $\varepsilon $ is a homomorphism. As in the
derivation of (7), we now have

\begin{eqnarray}
(a_{1},a_{2},a_{3}) &=&\left( a_{1},a_{2},\frac{3}{2}\frac{(a_{2})^{2}}{a_{1}%
}\right) \left( a_{1},a_{2},\frac{3}{2}\frac{(a_{2})^{2}}{a_{1}}\right)
^{-1}(a_{1},a_{2},a_{3})  \notag \\
&=&\left( a_{1},a_{2},\frac{3}{2}\frac{(a_{2})^{2}}{a_{1}}\right) \left(
1,0,a_{3}-\frac{3}{2}\frac{(a_{2})^{2}}{a_{1}}\right)
\end{eqnarray}%
and we observe that $a_{3}-\frac{3}{2}\frac{(a_{2})^{2}}{a_{1}}$ is the
defining formula for the Schwarzian derivative! Thus (51) defines the
semidirect product

\begin{equation}
G_{3}(1)=G_{2}(1)\ltimes K_{3,1}(1)
\end{equation}%
where $K_{3,1}(1)=\mathbb{F}$, but with a fundamental difference: (7) holds
for all $n\geq 1$ whereas (51) holds only for $n=1$, that is, we can define
an affine structure on any smooth manifold whereas to define the above
structure we must have, at least for the moment, $\dim M=1.$ We can now
define the groupoid $\mathcal{G}_{3}$, but what will we be dealing with? So
it seems reasonable to ask

$\mathbf{Q1:}$ Do all splittings inside (2) arise from Klein geometries?

On the infinesimal level, $\mathbf{Q1}$ asks whether all finite dimensional
Lie subalgebras of the Lie algebra of formal vector fields are induced by
Klein geometries $(\mathfrak{g},\mathfrak{h}).$ There is a striking relation
of (51) to Riccati equation ([1]).

Let $(G,M)$ be a Klein geometry such that $H\subset G$ is an algebraic
subgroup. We do not know whether the injection of $H$ inside $G_{m+1}(n)$
given by Proposition 5.5 in [3] imbeds $H$ as an algebraic subgroup. As
remarked in Section 2, we can always express $G_{m+1}(n)/H$ locally as a
zero set, but the problem is to find $F$ explicitly.

$\mathbf{Q2:}$ Given $(G,M),$ is there an effective algorithm for deriving
the defining equations of $\mathcal{G}_{m+1}$ in some canonical form ?

Recently the moving frame method of Cartan is perfected to its final form in
the remarkable papers [11], [12] in terms of a powerful and concrete
algorithm (see also [21] for a simplified overview). We believe that this
algorithm will play a decisive role in answering $\mathbf{Q2.}$ In fact, we
believe that the generalization of the Erlangen Program that we propose is a
special case and, we hope, also a geometric justification of the moving
frame method in [11], [12] in the case of transitive actions of Lie groups,
that is, pseudogroups of finite type.

Once we have the defining equations of $\mathcal{G}_{m+1}$, the rest is as
in this note, [2] and [25]: We linearize $\mathcal{G}_{m+1}$ to get the
algebroid $\mathfrak{G}_{m+1}.$ In these equations the top order one will
express $(m+1)^{th}$ order jet in terms of lower order jets but the lower
order equations may not be as innocent as in (28).

$\mathbf{Q3:}$ Is there an effective algorithm for computing the components
of the curvatures $\mathcal{R}$, $\mathfrak{R}$ ?

We again believe that the moving frame method will be of great help here.

To our embarassment, we are unable to prove the easy should be $\mathcal{R}%
=0\Rightarrow $ $\mathfrak{R}=0$ using (35), (38), (40), (42).

$\mathbf{Q4:}$ Give the invariant definitions of $\mathcal{R}$, $\mathfrak{R}
$ and prove the Lie's theorems $\mathcal{R}=0\Leftrightarrow $ $\mathfrak{R}%
=0$

The uniformization theorem in its full formulation gives rise to two
problems. First, when $\mathcal{R}=0$, the completeness of a pre-homogeneous
structure is defined in the same way. The equivalence $\mathcal{G}_{m+1}$
complete $\Leftrightarrow \mathfrak{G}_{m+1}$ complete is again not
difficult to show. However, metric completeness in $MRG$ suggests that
completeness should not need $\mathcal{R}=0.$ We believe that the
exponential map exists in some weak form without the assumption $\mathcal{R}%
=0$ (see remark 3) on pg. 22 of [2]) but are unable to make much progress
with it. Completeness is surely one of the most subtle concepts in the
theory.

$\mathbf{Q5:}$ What is the definition of completeness of a pre-homogeneous
structure in the presence of curvature?

Second, it is quite clear that the number of isomorphism classes of simply
connected Lie groups defining the vertex connection $\{G,N,H\}$ is equal to
the number of the possible uniformizing Klein geometries up to isomorphism.
We believe that this number and the vertex connection are uniquely
determined by the filtration (9) in [3]. For $m=0$, this number is infinite
since it is equal to the number of isomorphism classes of all Lie algebras.
We have seen that it is three for Riemannian structures and one for affine
structures.

$\mathbf{Q6:}$ Is this number always finite for $m\geq 1$ ?

The constant curvature condition (49), though it is quite natural and
geometric from the point of view of $MRG,$ is a total mystery for us and we
are unable to express the constant in (49) in terms of the structure
constants in (48).

$\mathbf{Q7:}$ \ Is (49) peculiar to Riemannian geometry or a particular
instance of a more general phenomenon ?

Finally we come to the surely most subtle part of the theory.

$\mathbf{Q8:}$ Develop the theory of characteristic classes (both primary
and secondary) for pre-homogeneous structures.

We believe that the variational (bi)complex ([33], [4], [31],[20]) and the
more recent invariant variational (bi)complex ([4], [5],[16]) will play a
fundamental role in such development.

\bigskip

\textbf{References}

\bigskip

[1] E.Abado\u{g}lu: 3-dimensional local Lie groups and bi-hamiltonian
systems, preprint, 2010

[2] E.Abado\u{g}lu, E.Orta\c{c}gil: Intrinsic characteristic classes of a
local Lie group, to appear in Portugaliae Mathematica

[3] E.Abado\u{g}lu, E.Orta\c{c}gil, F.\"{O}zt\"{u}rk: Klein geometries,
parabolic geometries and differential equations of finite type, J. Lie
Theory 18, (2008) no.1, 67-82

[4] I.M.Anderson: Introduction to the variational bicomplex.\ mathematical
aspects of classical field theory (Seattle, WA, 1991) 51-73, Contemp. Math.,
132, Amer. Math. Soc., Providence, RI

[5] I.M.Anderson, J. Pohjanpelto: The cohomology of invariant variational
bicomplexes, Geometric and algebraic structures in differential equations.
Acta Appl. Math. 41 (1995), no. 1-3, 3-19

[6] A.D.Blaom: Geometric structures as deformed infinitesimal symmetries,
Trans. Amer. Soc. 358 (2006) no.8, 3651-3671 (electronic)

[7] R.L.Bryant, S.S.Chern, R.B.Gardner, H.L.Goldschmidt, P.A.Griffiths:
Exterior differential systems, Mathematical Sciences Research Institute
Publications, !8. Springer-Verlag, New York, 1991

[8] A.Cap, J.Slovak, V.Soucek: Bernstein-Gelfand-Gelfand sequences, Ann. of
Math. 154 (2001), 97-113

[9] L.P.Eisenhart: Continuous groups of transformations, Dover Publications,
Inc., New York 1961

[10] L.P.Eisenhart: Riemannian geometry, 2d printing, Princeton University
Press, Princeton, N.J., 1949

[11] M.Fels, P.J.Olver: Moving coframes. I. A practical algorithm. Acta
Appl. Math. 51 (1998), no. 2, 161-213

[12] M.Fels, P.J.Olver: Moving coframes. II. Regularization and theoretical
foundations. Acta Appl. Math. 55 (1999), no.2, 127-208

[13] V.Guillemin: The integrability problem for G-structures, Trans. A.M.S.
116 (1965), 544-560

[14] S.Kobayashi, K.Nomizu: Foundations of differential geometry, Vol.1,
Interscience Publishers, John Wiley \& Sons, New York- London

[15] S.Kobayashi: Transformation groups in differential geometry, Ergebnisse
der Mathematik und ihrer Grenzgebiete, Band 70. Springer-Verlag, New
York-Heidelberg, 1972

[16] I.Kogan, P.J.Olver: The invariant variational bicomplex. The geometric
study of differential equations (Washington DC, 2000), 131-144, Contemp.
Math., 285, Amer. Math. Soc., Providence, RI, 2001

[17] D.Laugwitz: Differential and Riemannian geometry. Translated by Fritz
Steinhardt, Academic Press, New York-London, 1965

[18] K.MacKenzie: Lie groupoids and Lie algebroids in differential geometry,
London Mathematical Society Lecture Note Series, 124, Cambridge University
Press, Cambridge, 1987

[19] P.J.Olver: Nonassociative local Lie groups, J. Lie Theory 6 (1996) 23-51

[20] P.J.Olver: Applications of Lie groups to differential equations,
Graduate Texts in Mathematics, Springer-Verlag, New York, 1986

[21] P.J.Olver: An introduction to moving frames. Geometry, integrability
and quantization, 67-80, Softex, Sofia, 2004

[22] P.J.Olver: Symmetry groups and group invariant solutions of partial
differential equations, J. Differential Geom. 14 (1979), no. 4, 497-542

[23] P.J.Olver, J.Pohjanpelto: Maurer-Cartan forms and the structure of Lie
pseudo-groups. Selecta Math. (N.S.) 11 (2005), no. 1, 99-126

[24] E.Orta\c{c}gil: The heritage of Sophus Lie and Felix Klein: Geometry
via transformation groups, arXiv.org, math.DG/o6o4223, 2006

[25] J.F.Pommaret: Systems of partial differential equations and Lie
pseudogroups, with a preface by Andre Lichnerowicz, Mathematics and its
Applications, 14, Gordon \& Breach, Science Publishers, New York, 1978

[26] J.F.Pommaret: Partial differential equations and group theory, New
Perspectives for Applications,, 293, Kluwer Academic Publishers, 1994

[27] B.L.Reinhart: Some remarks on the structure of the Lie algebra of
formal vector fields, Transversal structure of foliations (Tolouse, 1982),
Asterisque, No. 116 (1984), 190-194

[28] A.Kumpera, D.C.Spencer: Lie equations, Vol.1: General theory, Annals of
Mathematics Studies, No.73, Princeton University Press, Princeton, NJ,
University of Tokyo Press, 1972

[29] R.Sharpe: Differential geometry, Cartan's generalization of Klein's
Erlangen program, with a foeward by S.S.Chern, Graduate Texts in
Mathematics, 166, Springer-Verlag, New York, 1997

[30] C.L.Terng: Natural vector bundles and natural differential operators,
Amer. J. Math., 100 (1978), no. 4, 775-828

[31] T.Tsujishita: On the variational bicomplexes associated to differential
equations, Osaka J. Math. 19 (1982), 311-363

[32] O.Veblen: Invariants of quadratic differential forms, Cambridge
University Press, 1962

[33] A.M.Vinogradov: A spectral sequence that is connected with a nonlinear
differential equation, and the algebraic-geometry foundations of the
Lagrange field theory with constraint (Russian) Dokl Akad. Nauk. SSSR 238
(1978), no.5, 1028-1031

[34] K.Yano: The theory of Lie derivatives and its applications,
North-Holland Publishing Co., Amsterdam, Noordhoff Ltd., Groningen,
Interscience Publishers Inc., New York, 1957

\bigskip

Erc\"{u}ment Orta\c{c}gil, Mathematics Department, Bo\u{g}azi\c{c}i
University, Bebek, 34342, Istanbul, Turkey

e-mail: ortacgil(@)boun.edu.tr

\end{document}